\documentclass[a4paper,pdftex,reqno,10pt]{amsart}

\allowdisplaybreaks[1]

\usepackage{geometry}

\usepackage{graphicx}
\graphicspath{{Figures/}}

\usepackage{times}
\usepackage{courier}
\usepackage[euler-digits]{eulervm}
\usepackage{eurosym}
\usepackage{hyperref}
\usepackage{amssymb,amsmath}

\usepackage{epsfig, graphicx}
\theoremstyle{plain}
\newtheorem{Theorem}{Theorem}[section]

\newtheorem{Conjecture}[Theorem]{Conjecture}

\newtheorem{Lemma}[Theorem]{Lemma}
\newtheorem{Construction}[Theorem]{Construction}

\newtheorem{Algorithm}[Theorem]{Algorithm}
\newenvironment{Proof}
{\begin{trivlist}\item[]{{\sc Proof.}}}{\hfill{$\square$}\noindent\end{trivlist}}

\theoremstyle{definition}
\newtheorem{Definition}[Theorem]{Definition}

\theoremstyle{remark}

\newcommand{\Tsm}{\hspace*{0.6cm}}

\begin{document}


\title[Maximal integral point sets over $\mathbf{\mathbb{Z}^2}$]{Maximal integral point sets over $\mathbf{\mathbb{Z}^2}$}

\author{Andrey Radoslavov Antonov}
\address{Andrey Radoslavov Antonov\\Department of Mathematics \\University of Chemical Technology and Metallurgy - Sofia \\Bulgaria}
\email{andrio@uctm.edu}

\author{Sascha Kurz}

\address{Sascha Kurz\\Department of Mathematics, Physic and Informatics\\University of Bayreuth\\Germany}
\email{sascha.kurz@uni-bayreuth.de}

\begin{abstract}
  Geometrical objects with integral side lengths have fascinated mathematicians through the ages. We call a set
  $P=\{p_1,\dots,p_n\}\subset\mathbb{Z}^2$ a maximal integral point set over $\mathbb{Z}^2$ if all pairwise
  distances are integral and every additional point $p_{n+1}$ destroys this property. Here we consider such sets
  for a given cardinality and with minimum possible diameter. We determine some exact values via exhaustive search
  and give several constructions for arbitrary cardinalities. Since we cannot guarantee the maximality in these cases
  we describe an algorithm to prove or disprove the maximality of a given integral point set. We additionally consider
  restrictions as no three points on a line and no four points on a circle.
\end{abstract}


\maketitle

\section{Introduction}

Geometrical objects with integral side lengths have fascinated mathematicians through the ages. A very early example is the Pythagorean triangle with side lengths $3$, $4$, and $5$. A universal framework for
most of these objects are integral point sets. By an integral point set we understand a set of $n$ points in an $m$ dimensional Euclidean vector space $\mathbb{E}^m$, where the pairwise distances between the points are integral. Those integral point sets
were studied by many authors, see \cite{integral_distances_in_point_sets} for an overview. From a combinatorial point of view
for a given cardinality $n$ and a given dimension $m$ the question on the minimum possible diameter $d(n,m)$, this is the largest
distance between any two points, arises, see \cite{phd_kurz,paper_laue,paper_alfred} for an overview.

To obtain some interesting discrete structures one could also require some additional properties. One possibility is to request,
that besides the distances also the coordinates must be integral. Another classical possibility is to forbid subsets of three points on a line or four points on a circle. The question of P. Erd\H{o}s whether there exists a set of seven points in the plane with no
three points on a line, no four points on a circle, and pairwise integral distances, has recently been answered positively, see \cite{paper_kreisel}. If all three mentioned additional properties are required simultaneously one speaks of $n_m$-clusters, see \cite{0676.52006}. In this article we request that besides the distances also the coordinates of the point sets are integral and restrict ourselves to dimension $2$. Additionally we consider the cases where no three points are on a line or no four points are on a circle.

In finite geometry one is sometimes interested in point configurations which are maximal with respect to some property. This means that it is not possible to add a point without destroying the requested property. Here we consider integral point sets which are maximal, meaning that there does not exist an additional point $x$ with integral distances to the other points of the point set. 

\subsection{Related work}
There have been extensive studies on integral point sets in Euclidean spaces. Some authors also consider other spaces, e.~g. Banach spaces \cite{banach}, integral point sets over rings \cite{paper_axel_2}, or integral point sets over finite fields \cite{algo,inclusion_maximal,integral_over_finite_fields}. In \cite{gauss_integers} the authors consider integral point sets over $\mathbb{Z}^2$ and conjecture some examples to be maximal. As an answer to their open problems in \cite{paper_axel} the authors describe an algorithm to prove the maximality of a given integral point set and prove the conjectures of \cite{gauss_integers}.

\subsection{Our contribution}
In this paper we describe algorithms to efficiently test integral point sets for maximality and to determine possible extension points. To deal with the isomorphism problem we describe an algorithm which transforms a given plane integral point set into a normal form in $O\left(n^2\right)$ time, where $n$ is the cardinality of the point set. We give several constructions of integral point sets over $\mathbb{Z}^2$ which have a given cardinality and fulfill additional conditions such as that there are {\lq\lq}no three points on a line{\rq\rq} or {\lq\lq}no four points on a circle{\rq\rq}. Although we cannot prove the maximality of the point sets obtained with the proposed constructions in general, we conjecture this property for many of our constructions. By exhaustive search we have determined some exact minimum diameters of integral point sets over $\mathbb{Z}^2$ with given cardinality and with or without additional conditions. We give constructive upper bounds in most cases and conjecture them to be the exact values.

\subsection{Outline of the paper}
In Section \ref{sec_exhaustive_generation} we state the basic definitions and in Section \ref{sec_basics} we describe the basic algorithms to deal with maximal integral point sets over $\mathbb{Z}^2$. These include an algorithm to exhaustively generate Heronian triangles up to isomorphism, an algorithm to determine all possible embeddings of an Heronian triangle on the integer grid $\mathbb{Z}^2$, and an algorithm that determines all points of $\mathbb{Z}^2$ which have integral distances to three given points in $\mathbb{Z}^2$ with pairwise integral distances. The last mentioned algorithm enables us to algorithmically prove or disprove the maximality of a given integral point set. Since we intend to consider integral point sets up to isomorphism we introduce normal forms of integral point sets and algorithms to obtain them in Section \ref{sec_normal_forms}. We deal with the key question of maximal integral point sets over $\mathbb{Z}^2$ with given cardinality and minimum diameter in Section \ref{sec_maximal_integral_point_sets}.
Several constructions for maximal integral point sets, where the maximality is not guaranteed but very likely, are described in Section \ref{sec_constructions}. In Section \ref{sec_further_conditions} we deal with additional properties as {\lq\lq}no three points on a line{\rq\rq} and {\lq\lq}no four points on a circle{\rq\rq}. We finish with a short conclusion and an outlook in Section \ref{sec_outlook}.

\section{Basics}
\label{sec_basics}

\noindent
\begin{Definition}
  An integral point set over $\mathbb{Z}^2$ is a non-collinear set $\mathcal{P}$ of $n$ points in the integer grid $\mathbb{Z}^2$,
  where the points have pairwise integral distances.
\end{Definition}

For brevity we only speak of integral point sets and assume that the coordinates of the points are integral numbers, too.

\begin{Definition}
  We call an integral point set $\mathcal{P}$ over $\mathbb{Z}^2$ maximal if for every $x\in\mathbb{Z}^2\backslash\mathcal{P}$ the
  point set $\mathcal{P}\cup \{x\}$ is not an integral point set.
\end{Definition}

The existence of maximal integral point sets in the plane is guaranteed  by a famous theorem of N.H..~Anning and P.~Erd\H{o}s, respectively its proof.

\begin{Theorem}
  An infinite set $\mathcal{P}$ of points in the Euclidean space $\mathbb{E}^m$ with pairwise integral distances is situated
  on a line. \cite{ErdoesAnning1,ErdoesAnning2}
\end{Theorem}
\begin{Proof}
  We only prove the statement for dimension $m=2$, as the generalization is obvious. If $A$, $B$, and $C$ are three points not on a
  line, we set $k=\max\left\{\overline{AC},\overline{BC}\right\}$ and consider points $P$ such that $|\overline{PA}-\overline{PC}|$
  and
  $|\overline{PB}-\overline{PC}|$ are integral. Due to the triangle inequalities the attained values are in $\{0,1,\dots,k\}$. Thus
  the point $P$ lies on the intersection of two distinct hyperbolas, where we have at most $k+1$ choices for each hyperbola. Thus
  there are at most $4(k+1)^2$ possible locations for the point $P$.
\end{Proof}

\begin{figure}[htp]
  \begin{center}
    \setlength{\unitlength}{0.50cm}
    \begin{picture}(6,10)
      \put(3,2){\circle*{1}}
      \put(3.4,2){$(0,-4)$}
      \put(0,6){\circle*{1}}
      \put(-1.8,6.55){$(-3,0)$}
      \put(3,6){\circle*{1}}
      \put(3.4,6.1){$(0,0)$}
      \put(6,6){\circle*{1}}
      \put(6,6.55){$(3,0)$}
      \put(3,10){\circle*{1}}
      \put(3.4,10){$(0,4)$}
      \put(0,6){\line(1,0){6}}
      \put(3,2){\line(0,1){8}}
      \put(0,6){\line(3,4){3}}
      \put(0,6){\line(3,-4){3}}
      \put(6,6){\line(-3,4){3}}
      \put(6,6){\line(-3,-4){3}}
    \end{picture}
    \quad\quad\quad
    \setlength{\unitlength}{0.25cm}
    \begin{picture}(25,24)
      \put(0,12){\circle*{2}}
      \put(-2.9,9.8){$(0,12)$}
      \put(9,0){\circle*{2}}
      \put(5.0,-0.2){$(9,0)$}
      \put(16,0){\circle*{2}}
      \put(16.9,-0.2){$(16,0)$}
      \put(9,24){\circle*{2}}
      \put(4.2,23.8){$(9,24)$}
      \put(16,24){\circle*{2}}
      \put(16.9,23.8){$(16,24)$}
      \put(25,12){\circle*{2}}
      \put(25.9,11.7){$(25,12)$}
      \put(0,12){\line(1,0){25}}
      \put(0,12){\line(3,4){9}}
      \put(0,12){\line(3,-4){9}}
      \put(25,12){\line(-3,4){9}}
      \put(25,12){\line(-3,-4){9}}
      \put(9,0){\line(1,0){7}}
      \put(9,24){\line(1,0){7}}
      \put(9,0){\line(0,1){24}}
      \put(16,0){\line(0,1){24}}
      \put(0,12){\line(4,3){16}}
      \put(0,12){\line(4,-3){16}}
      \put(25,12){\line(-4,3){16}}
      \put(25,12){\line(-4,-3){16}}
      \qbezier(9,0)(12.5,12)(16,24)
      \qbezier(9,24)(12.5,12)(16,0)
    \end{picture}
  \end{center}
  \caption{Examples of maximal integral point sets.}
  \label{fig_examples_old}
\end{figure}
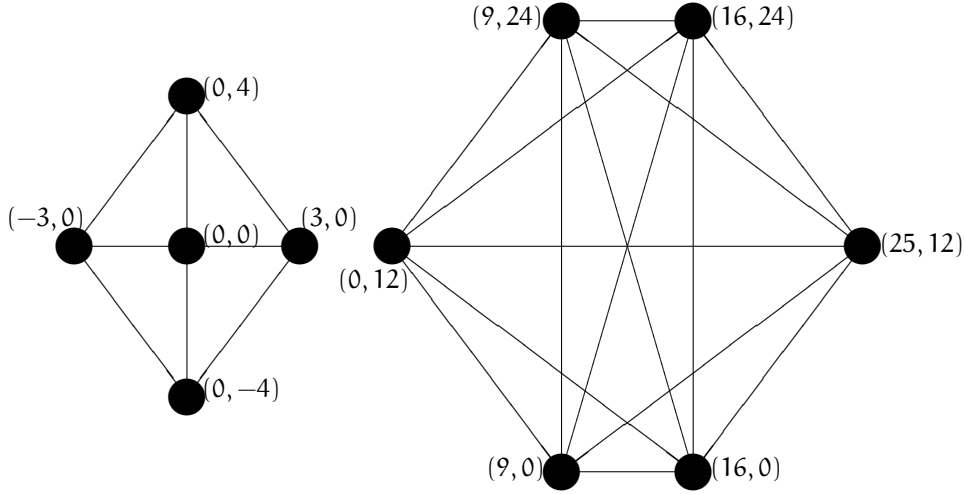

This proof can clearly be converted into a constructive algorithm. Given three points $A=(x_A,y_A)$, $B=(x_B,y_B)$, and $C=(x_C,y_C)$ in $\mathcal{P}\subset\mathbb{Z}^2$, which are not on a line, the problem of determing points $P=(x_P,y_P)$ at integral distance to $A$, $B$, and $C$ is reduced to the problem of solving the equation system
\begin{equation}
  \label{eq_fourth_point}
  \left|
    \begin{array}{rcl}
      \sqrt{(x_A-x_P)^2+(y_A-y_P)^2}-\sqrt{(x_C-x_P)^2+(y_C-y_P)^2}&=&d_1\\
      \sqrt{(x_B-x_P)^2+(y_B-y_P)^2}-\sqrt{(x_C-x_P)^2+(y_C-y_P)^2}&=&d_2\\
    \end{array}
  \right|,
\end{equation}
where $d_1\in\left\{-\overline{AC},\dots,\overline{AC}\right\}\subset\mathbb{Z}$ and $d_2\in\left\{-\overline{BC},\dots,\overline{BC}\right\}\subset\mathbb{Z}$. If there exists no integral solution in $\mathbb{Z}^2\backslash\mathcal{P}$, then the point set $\mathcal{P}$ is maximal. This algorithm was already used in \cite{paper_axel} to prove the maximality of the integral point sets of Figure \ref{fig_examples_old}.

Since this algorithm is essential for our article we will go into the details how to solve equation system \ref{eq_fourth_point}. To get rid of some of the square roots we add $\sqrt{(x_C-x_P)^2+(y_C-y_P)^2}$ on both sides and square the expressions afterwards:
$$
  \left|
    \begin{array}{rcl}
      (x_A-x_P)^2+(y_A-y_P)^2 &=& d_1^2+2d_1\sqrt{(x_C-x_P)^2+(y_C-y_P)^2}+(x_C-x_P)^2+(y_C-y_P)^2\\
      (x_B-x_P)^2+(y_B-y_P)^2 &=& d_2^2+2d_2\sqrt{(x_C-x_P)^2+(y_C-y_P)^2}+(x_C-x_P)^2+(y_C-y_P)^2\\
    \end{array}
  \right|.
$$
Rearranging yields
\begin{equation}
  \label{eq_last_step}
  \left|
    \begin{array}{rcl}
      (x_A^2+y_A^2-x_C^2-y_C^2-d_1^2)+2(x_C-x_A)x_P+2(y_C-y_A)y_P &=&
      2d_1\sqrt{(x_C-x_P)^2+(y_C-y_P)^2}\\
      (x_B^2+y_B^2-x_C^2-y_C^2-d_2^2)+2(x_C-x_B)x_P+2(y_C-y_B)y_P  &=&
      2d_2\sqrt{(x_C-x_P)^2+(y_C-y_P)^2}\\
    \end{array}
  \right|.
\end{equation}
If $d_1=0$ then the first equation corresponds to a linear equation 
\begin{equation}
  \label{eq_lin}
  c_1x_P+c_2y_P+c_3=0,
\end{equation}
where not both $c_1$ and $c_2$ are equal to zero, since $A\neq C$. If we square the second equation of (\ref{eq_last_step}) we can substitute one variable using equation (\ref{eq_lin}) and obtain a quadratic equation in one variable, which can be easily solved. The case, where $d_2=0$ is similar. Here we use the second equation of (\ref{eq_last_step}) to obtain equation (\ref{eq_lin}) (we have $c_1\neq 0$ or $c_2\neq 0$ due to $B\neq C$), and substitute it into the squared version of the first equation to obtain the quadratic equation in one variable. In the remaining case we have $d_1,d_2\neq 0$. Here we subtract $d_1$ times the second equation of (\ref{eq_last_step}) from $d_2$ times the first equation of (\ref{eq_last_step}) to obtain equation (\ref{eq_lin}) (we have $c_1\neq 0$ or $c_2\neq 0$ since the points $A$, $B$, and $C$ are not located on a line). Now we can square one of the two equations of (\ref{eq_last_step}) and subsitute one variables using equation (\ref{eq_lin}). Again we end up with a quadratic equation in one variable. At the end we have to check if the obtained values $(x_P,y_P)$ are solutions of the original equation system (\ref{eq_fourth_point}).

\begin{Definition}
  For an integral point set $\mathcal{P}$ its diameter $\mbox{diam}(\mathcal{P})$ is given by the largest distance between
  a pair of its points.
\end{Definition}

We remark that the left integral point set  of Figure \ref{fig_examples_old} has diameter $8$ and the right integral point set of Figure \ref{fig_examples_old} has diameter $25$.

\section{Exhaustive generation of maximal integral point sets}
\label{sec_exhaustive_generation}

\noindent
To obtain interresting examples of maximal integral point sets we utilize computers to exhaustively generate maximal integral point sets. In the following we will describe the algorithm used. For a given diameter $d$ we loop over all non-isomorphic Heronian triangles (having integral side lengths and integral area) $\Delta=(a,b,c)$ with diameter $d=\max\{a,b,c\}$. Utilizing the Heron formula
\begin{equation}
  A=\frac{\sqrt{(a+b+c)(a+b-c)(a-b+c)(-a+b+c)}}{4}
\end{equation}
for the area of a triangle we can generate this list e.g. by the following short algorithm:
\begin{Algorithm}{(Generation of Heronian triangles)}\\
  \label{algo_triangles}
  \noindent
  \textit{input:} diameter $d$\\
  \textit{output:} complete list of Heronian triangles with diameter $d$ up to isomorphism\\
  {\bf begin}\\
  \Tsm $a=d$\\
  \Tsm {\bf for} $b=\left\lfloor\frac{a+2}{2}\right\rfloor,\dots,a$ {\bf do}\\
  \Tsm\Tsm {\bf for} $c=a+1-b,\dots,b$ {\bf do}\\
  \Tsm\Tsm\Tsm {\bf if} $\frac{\sqrt{(a+b+c)(a+b-c)(a-b+c)(-a+b+c)}}{4}\in\mathbb{Z}$ {\bf then}\\
  \Tsm\Tsm\Tsm\Tsm output $(a,b,c)$\\ 
  {\bf end}
\end{Algorithm}

For a more sophisticated and efficient algorithm we refer to \cite{herontriangles}. The next step is to embed a given Heronian triangle $\Delta=(a,b,c)$ in the plane integer grid $\mathbb{Z}^2$. Here we can utilize two conjectures, which are theorems for dimension $m=2$, see e.g. \cite{fricke}.

\begin{Conjecture}
  \label{conj_embedding}
  Let $\mathcal{P}\subset\mathbb{Q}^m$ be a finite set of points such that the distances between any two points of $\mathcal{P}$
  are integers. In this case one can find an Euclidean motion $T$ such that $T(\mathcal{P})\subset\mathcal{P}^m$. 
\end{Conjecture}

\begin{Conjecture}
  Let $\mathcal{P}\subset\mathbb{Z}^m$ be a finite set of points such that the distances between any two points of $\mathcal{P}$
  are integers and divisible by an integer $k$. In this case one can find a set $\mathcal{P}'\subset\mathbb{Z}^m$ such that
  $\mathcal{P}'\cdot k$ (the set $\mathcal{P}'$ scaled by a factor $k$) is congruent to $\mathcal{P}$.
\end{Conjecture}

Since Conjecture \ref{conj_embedding} is a well known theorem for dimension $m=2$, see e.g. \cite{fricke}, for every Heronian triangle $\Delta(a,b,c)$ there exists an embedding in the plane integer grid $\mathbb{Z}^2$. We remark that there may be several embeddings for the same triangle $\Delta=(a,b,c)$, which lead to different results. If we consider the number of points $(x_P,y_P)\in\mathbb{Z}^2\backslash\mathcal{E}$ which are at integral distance to an embedded triangle $\mathcal{E}=\{(x_A,y_A),(x_B,y_B),(x_C,y_C)\}$, we can distinguish three different embeddings of the Heronian triangle $\Delta_1=(25,20,15)$. The embedding $\mathcal{E}_1=\{(0,0),(0,25),(12,16)\}$ of $\Delta_1$ yields $12$ points $(x_P,y_P)$ at integral distance to the corners of $\Delta_1$ given by $\mathcal{E}_1$. For the embedding $\mathcal{E}_2=\{(0,0),(15,20),(0,20)\}$
we obtain $16$ such points, and for the embedding $\mathcal{E}_3=\{(0,0),(7,24),(16,12)\}$ we obtain only $5$ such points. Determining the possible embeddings of a given Heronian triangle $\Delta=(a,b,c)$ is a rather easy task. W.l.o.g. we assume $a=\max\{a,b,c\}$ and $x_B=0=y_B$. Since the point $(x_C,y_C)$ is at distance $a$ to the point $(x_B,y_B)$, we have to solve the Diophantine equation
$$
  x_C^2+y_C^2=a^2
$$
in integers. This is a well known problem. One might even store for each small number (e.g. $a\le 10\,000$) $a\in\mathbb{N}$ a list of the corresponding solutions $(x_C,y_C)$. Now the coordinates of the remaining point $A$ are given as solutions of the equation system
\begin{equation}
  \label{eq_third_point}
  \left|
    \begin{array}{rcl}
      (x_B-x_A)^2+(y_B-y_A)^2 &=& c^2 \\
      (x_C-x_A)^2+(y_C-y_A)^2 &=& b^2
    \end{array}
  \right|,
\end{equation}
which can be easily solved. As an algorithm for the embedding of an Heronian triangle in $\mathbb{Z}^2$ we obtain:
\begin{Algorithm}{(Embedding of an Heronian Triangle)}\\
  \label{algo_embedding}
  \noindent
  \textit{input:} Heronian Triangle $\Delta=(a,b,c)$\\
  \textit{output:} complete list of different embeddings of $\Delta$ in $\mathbb{Z}^2$\\
  {\bf begin}\\
  \Tsm $x_B=0$,\,$y_B=0$\\
  \Tsm {\bf loop over} the integer solutions $(x_C,y_C)$ of $x_C^2+y_C^2=a^2$ {\bf do}\\
  \Tsm\Tsm {\bf loop over} the integer solutions $(x_A,y_A)$ of equation system (\ref{eq_third_point}) {\bf do}\\
  \Tsm\Tsm\Tsm output $\{(x_A,y_A),(x_B,y_B),(x_C,y_C)\}$\\
  {\bf end}
\end{Algorithm}
The next step is to determine the points $(x_P,y_P)\in\mathbb{Z}^2$ which are at integral distance to a given embedded triangle  $\{(x_A,y_A),(x_B,y_B),(x_C,y_C)\}$:
\begin{Algorithm}{(Enlargement of an embedded triangle)}\\
  \label{algo_enlargement}
  \noindent
  \textit{input:} Embedded triangle $\mathcal{E}=\{(x_A,y_A),(x_B,y_B),(x_C,y_C)\}\subset\mathbb{Z}^2$\\
  \textit{output:} complete list of points $(x_P,y_P)\in\mathbb{Z}^2\backslash\mathcal{E}$ which are at
                   integral distance to $\mathcal{E}$\\
  {\bf begin}\\
  \Tsm {\bf loop over} the integer solutions $(x_P,y_P)$ of equation system (\ref{eq_fourth_point}) {\bf do}\\
  \Tsm\Tsm {\bf if} $(x_P,y_P)\notin\mathcal{E}$ {\bf then}\\
  \Tsm\Tsm\Tsm output $(x_P,y_P)$\\
  {\bf end}
\end{Algorithm}
We remark that the previous algorithms have to be implemented using an arithmetic which is able to do integer calculations with unlimited precision, since the occurring numbers can increase very quickly. We have utilized the software package \textsc{CLN} \cite{cln} for this purpose.

Now we utilize the set of points given by Algorithm \ref{algo_enlargement} to build up a graph $\mathcal{G}(\mathcal{E})$. The vertices are
given by the possible points $(x_P,y_P)$. Two points $\left(x_{P_1},y_{P_1}\right)$ and $\left(x_{P_2},y_{P_2}\right)$ are connected
by an edge if and only if $\sqrt{(x_{P_1}-x_{P_2})^2+(y_{P_1}-y_{P_2})^2}$ is a positive integer. A complete subgraph of $\mathcal{G}(\mathcal{E})$ is called a clique. A clique $\mathcal{C}_1$ is called maximal if it is not properly contained in another clique $\mathcal{C}_2$ of $\mathcal{G}(\mathcal{E})$. Clearly the cliques of $\mathcal{G}(\mathcal{E})$ are in bijection to integral point sets $\mathcal{P}\subset\mathbb{Z}^2$ containing $\mathcal{E}$ as a subset. The same statement holds for maximal cliques of $\mathcal{G}(\mathcal{E})$ and maximal integral point sets $\mathcal{P}\subset\mathbb{Z}^2$ containing $\mathcal{E}$ as a subset. Thus we can use a clique-search package as \textsc{Cliquer} \cite{cliquer} to exhaustively generate maximal integral point sets $\mathcal{M}$ over $\mathbb{Z}^2$.

Let us consider an example. If we apply our algorithm on the embedded triangle $$\mathcal{E}_2=\{(0,0),(15,20), (0,20)\}$$ with diameter $25$, we obtain a set
\begin{eqnarray*}
  \{(0,28),(0,40),(0,56),(0,132),(0,-92),(0,-16),(0,12),(-15,20),\\(15,0),(-21,20),(105,-36),(21,20),(-48,20),(48,20),(-99,20)\}
\end{eqnarray*}
of $16$ possible points to enlarge the integral point set $\mathcal{E}_2$. The clique-search program \textsc{Cliquer} determines five maximal cliques which correspond to the following five maximal integral point sets:
\begin{eqnarray*}
  \mathcal{M}_1 &=& \{(0,0),(15,20),(0,20),(15,0)\},\\
  \mathcal{M}_2 &=& \{(0,0),(15,20),(0,20),(0,-92),(105,-36)\},\\
  \mathcal{M}_3 &=& \{(0,0),(15,20),(0,20),(0,40),(0,56),(0,-16),(-15,20),(-48,20),(48,20)\},\\
  \mathcal{M}_4 &=& \{(0,0),(15,20),(0,20),(0,40),(-15,20),(-21,20),(21,20),(-48,20),(48,20),\\
                & & (-99,20),(99,20)\},\text{ and}\\
  \mathcal{M}_5 &=& \{(0,0),(15,20),(0,20),(0,28),(0,40),(0,56),(0,132),(0,-92),(0,-16),(0,12),\\
                & & (-15,20)\}.
\end{eqnarray*}
It is interesting to have a look at the cardinalities and diameters of these maximal integral point sets. We have $|\mathcal{M}_1|=4$, $\mbox{diam}(\mathcal{M}_1)=25$, $|\mathcal{M}_2|=5$, $\mbox{diam}(\mathcal{M}_2)=119$, $|\mathcal{M}_3|=9$, $\mbox{diam}(\mathcal{M}_3)=96$, $|\mathcal{M}_4|=11$, $\mbox{diam}(\mathcal{M}_1)=198$, $|\mathcal{M}_5|=11$, and $\mbox{diam}(\mathcal{M}_5)=224$. Although we start with a point set $\mathcal{E}_2$ of small diameter, the resulting maximal integral point sets $\mathcal{M}_i$ may have a large diameter. We are not aware of a formula to bound $\mbox{diam}(\mathcal{M})$ with respect to $\mbox{diam}(\mathcal{E})$. A second somewhat disappointing fact of our algorithm is, that each subset $\mathcal{E}'$ of
three non-collinear points of an maximal integral point set $\mathcal{M}$ produces $\mathcal{M}$. Thus our algorithm produces many identical copies of maximal integral point sets with large cardinality. We will deal with this fact and the isomorphism problem in the next section.

The algorithms described in this section focus on the maximality of the integral point sets. They should not be used to exhaustively generate all maximal integral point sets up to a given diameter. To perform this task the algorithms to exhaustively generate integral point sets with or without additional properties are better suited, see \cite{phd_kurz,paper_alfred}, and ignore the maximality condition in the first run. All integral point sets with required cardinalities and small diameters can then be tested if they are maximal.

\section{Normal forms and automorphisms for integral point sets over $\mathbf{\mathbb{Z}^2}$}
\label{sec_normal_forms}

\noindent
In this section we aim to consider isomorphisms which preserve certain properties of maximal integral point sets. Since a main property of an integral point set is the set of distances between its points we only consider distance-preserving isomorphisms, so called isometries. In the Euclidean plane the isometries are given by compositions of translations $T_{u,v}:\begin{pmatrix}x\\y\end{pmatrix}\mapsto \begin{pmatrix}x\\y\end{pmatrix}+\begin{pmatrix}u\\v\end{pmatrix}$, rotations $R_\theta:\begin{pmatrix}x\\y\end{pmatrix}\mapsto \begin{pmatrix}\cos\theta&-\sin\theta\\\sin\theta&\cos\theta\end{pmatrix}\cdot
\begin{pmatrix}x\\y\end{pmatrix}$, and reflections at one of the two axes. Each isometry can be written as $I_{t,O}:x\mapsto t+O\cdot x$, where $t\in\mathbb{R}^2$ is a translation vector and $O\in\mathbb{R}^{2\times 2}$ an orthogonal matrix. Next we restrict ourselves to mappings which map integral coordinates onto integral coordinates. Thus we have $t\in\mathbb{Z}^2$ and $O\in\mathbb{Z}^{2\times 2}$. Each such isometry $I_{t,O}$ maps integral point sets onto integral point sets. It is easy to figure out that there are only $8$ orthogonal matrices in $\mathbb{Z}^{2\times 2}$. So we define
$$
  \mbox{Aut}:=\left\{I_{t,O}\,:\, t\in\mathbb{Z}^2,\,O\in\left\{
  \begin{pmatrix}\pm 1&0\\0&\pm 1\end{pmatrix},\begin{pmatrix}\pm 1&0\\0&\mp 1\end{pmatrix},
  \begin{pmatrix}0&\pm 1\\\pm 1&0\end{pmatrix},\begin{pmatrix}0&\pm 1\\\mp 1&0\end{pmatrix}
  \right\}\right\}
$$
as the automorphism group of plane integral point sets.

We call two integral point sets $\mathcal{P}$ and $\mathcal{P}'$ isomorphic, if there exists a mapping $I_{t,O}\in\mbox{Aut}$ such that $I_{t,O}(\mathcal{P})=\mathcal{P}'$. So our aim is to develop an algorithm which can check whether two given integral point sets are isomorphic. For this purpose we want to use the technique of normal forms of discrete objects. This means that we have a function $\tau$ which fulfills the following: If $\mathcal{O}$ is the orbit of an integral point set $\mathcal{P}$ under the group $\mbox{Aut}$ then $\tau(\mathcal{P})=\tau(\mathcal{P}')$ for each $\mathcal{P}'\in\mathcal{O}$. Additionally for two integral point sets of different orbits the function $\tau$ should have different images. Having such a function $\tau$ at hand we can easily decide whether two integral point sets $\mathcal{P}$ and $\mathcal{P}'$ are isomorphic, by checking whether $\tau(\mathcal{P})=\tau(\mathcal{P}')$ or not.

In order to describe such a function $\tau$ we need to define a total ordering $\preceq$ on $\mathbb{Z}^2$:
\begin{itemize}
 \item[(1)] if $|a|<|c|$, then we set $\begin{pmatrix}a\\b\end{pmatrix}\prec\begin{pmatrix}c\\d\end{pmatrix}$,
 \item[(2)] if $a>0$, then we set $\begin{pmatrix}-a\\b\end{pmatrix}\prec\begin{pmatrix}a\\d\end{pmatrix}$,
 \item[(3)] if $|b|<|d|$, then we set $\begin{pmatrix}a\\b\end{pmatrix}\prec\begin{pmatrix}a\\d\end{pmatrix}$, and
 \item[(4)] if $b>0$, then we set $\begin{pmatrix}a\\-b\end{pmatrix}\prec\begin{pmatrix}a\\b\end{pmatrix}$
\end{itemize}
for all $a,b,c,d\in\mathbb{Z}$. We set $\begin{pmatrix}a\\b\end{pmatrix}=\begin{pmatrix}c\\d\end{pmatrix}$ if and only if
we have $a=c$ and $b=d$. By $x_1 \preceq x_2$ we mean $x_1\prec x_2$ or $x_1=x_2$. One of the properties of this total ordering $\preceq$ is, that we have $\begin{pmatrix}0\\0\end{pmatrix}\preceq x$ for all $x\in\mathbb{Z}^2$, so $\begin{pmatrix}0\\0\end{pmatrix}\preceq x$ is the smallest element in $\mathbb{Z}^2$. Using $\prec$ we can bijectively identify an integral point set $\mathcal{P}$ with a list $\mathcal{L}(\mathcal{P})$ of its points, which is sorted in ascending order with respect to $\preceq$. Now we extend our total ordering $\preceq$ onto such lists by utilizing the lexicographic ordering. This allows us to define our normalization function by
$$
  \tau(\mathcal{P})=\min_{\preceq}\left\{\mathcal{L}(\sigma(\mathcal{P}))\,:\,\sigma\in\mbox{Aut}\right\}.
$$
To obtain a finite algorithm for the determination of $\tau(\mathcal{P})$ we use the fact, that for every point set $\mathcal{P}\neq\emptyset$ the \textit{minimum} list-representation $\mathcal{L}(\sigma(\mathcal{P}))$ starts with $\begin{pmatrix}0\\0\end{pmatrix}$:
\begin{Algorithm}{(Normalization of an integral point set)}\\
  \label{algo_normalization}
  \noindent
  \textit{input:} integral point set $\mathcal{P}=\{p_1,\dots,p_n\}$\\
  \textit{output:} minimum list representation $\tau(\mathcal{P})$\\
  {\bf begin}\\
  \Tsm $champion=\mathcal{L}(\mathcal{P})$\\
  \Tsm $M_1=\begin{pmatrix}1&0\\0&1\end{pmatrix}$, $M_1=\begin{pmatrix}1&0\\0&-1\end{pmatrix}$,
       $M_3=\begin{pmatrix}-1&0\\0&1\end{pmatrix}$, $M_1=\begin{pmatrix}-1&0\\0&-1\end{pmatrix}$\\
  \Tsm $M_5=\begin{pmatrix}0&1\\1&0\end{pmatrix}$, $M_6=\begin{pmatrix}0&1\\-1&0\end{pmatrix}$,
       $M_7=\begin{pmatrix}0&-1\\1&0\end{pmatrix}$, $M_8=\begin{pmatrix}0&-1\\-1&0\end{pmatrix}$\\
  \Tsm {\bf for} $i=1,\dots,n$ {\bf do}\\
  \Tsm\Tsm {\bf for} $j=1,\dots,8$ {\bf do}\\
  \Tsm\Tsm\Tsm $tmp=\mathcal{L}(M_j\cdot\{p_1-p_i,\dots,p_n-p_i\})$\\
  \Tsm\Tsm\Tsm {\bf if} $tmp\prec champion$ {\bf then}\\
  \Tsm\Tsm\Tsm\Tsm $champion=tmp$\\
  \Tsm {\bf return} $champion$\\
  {\bf end}
\end{Algorithm}
We remark that Algorithm \ref{algo_normalization} runs in $O\left(n^2\right)$ time. As an example we consider the two integral point sets from Figure \ref{fig_examples_old}. Their normal forms or minimum list representations are given by
$$
  \left[
  \begin{pmatrix}0\\0\end{pmatrix},
  \begin{pmatrix}0\\-3\end{pmatrix},
  \begin{pmatrix}0\\3\end{pmatrix},
  \begin{pmatrix}-4\\0\end{pmatrix},
  \begin{pmatrix}4\\0\end{pmatrix}
  \right]
$$
and
$$
  \left[
  \begin{pmatrix}0\\0\end{pmatrix},
  \begin{pmatrix}0\\-7\end{pmatrix},
  \begin{pmatrix}-12\\9\end{pmatrix},
  \begin{pmatrix}-12\\-16\end{pmatrix},
  \begin{pmatrix}-24\\0\end{pmatrix},
  \begin{pmatrix}-24\\-7\end{pmatrix}
  \right],
$$
respectively.

For a given integral point set $\mathcal{P}$ there may exist rotation matrices $M\in\mathbb{R}^{2\times 2}$, such that $M(\mathcal{P})$ has integral coordinates, which are different from the eight orthogonal matrices in $\mathbb{Z}^{2\times 2}$. But for these matrices there is no 	guarantee for a proper extension $\mathcal{E}\supset\mathcal{P}$, which is also an integral point set over $\mathbb{Z}^2$, such that $M(\mathcal{E})$ has integral coordinates. Examples are given by the sets $\mathcal{E}_1$, $\mathcal{E}_2$, $\mathcal{E}_3$ in Section \ref{sec_exhaustive_generation}. This means that for a given maximal integral point set $\mathcal{M}$ over $\mathbb{Z}^2$ there can exist an orthogonal matrix $M\in\mathbb{R}^{2\times 2}$, such that $M(\mathcal{M})$ is also an integral point set over $\mathbb{Z}^2$, but which is not maximal.

We may call a maximal integral point set $\mathcal{M}$ over $\mathbb{Z}^2$ strongly maximal, if such a matrix $M$ does not exist. To check whether a given integral point set $\mathcal{P}$ is strongly maximal, we only have to consider all possible embeddings of $\mathcal{P}$ in $\mathbb{Z}^2$, which are finitely many. Another possibility is to slightly alter Algorithm \ref{algo_enlargement} by looping over the rational (instead of integral) solutions $(x_P,y_P)$ of equation system (\ref{eq_fourth_point}). Now the algorithm leads to point sets with integral distances and rational coordinates. But due to Conjecture \ref{conj_embedding} (which is a theorem for dimension $m=2$), there exist embeddings with integral coordinates.

To clear the situation with integral and rational coordinates we will have to give some facts from the general theory of integral point sets (without integral coordinates). So, let $\mathcal{P}$ be a set of points in the $m$-dimensional Euclidean space $\mathbb{E}^m$ with pairwise integral distances. By $\mathcal{S}\subseteq\mathcal{P}$ we denote an integral simplex, which is a set of $m+1$ points, and by $mbox{vol}_m(\mathcal{S})$ we denote the $m$-dimensional volume spanned by the $m+1$ points. Since the pairwise distances are integral we can write $\mbox{vol}_m(\mathcal{S})=q\cdot k$ with $q\in\mathbb{Q}$ and $k$ being a square free integer. If $\mbox{vol}_m(\mathcal{S})\neq 0$ the square free integer $k$ is unique and we set $\mbox{char}(\mathcal{S})=k$, which we call the characteristic of $\mathcal{S}$. Using this notation we can cite two results from \cite{characteristic}:

\begin{Theorem}
  \label{thm_characteristic}
  In an $m$-dimensional integral point set $\mathcal{P}$ all simplices $\mathcal{S}=\{v_0,v_1,\dots,v_{m}\}$
  with $\mbox{vol}_{m}(\mathcal{S})\neq 0$ have the same characteristic $\mbox{char}(\mathcal{S})=k$.
\end{Theorem}

So we can speak of \textit{the} characteristic $\mbox{char}(\mathcal{P})$ of an integral point set $\mathcal{P}$.

\begin{Lemma}
  \label{lemma_coordinates}
  An integral $m$-dimensional simplex $\mathcal{S}=\{v_0',v_1',\dots,v_{m}'\}$ with 
  distance matrix $D=(d_{i,j})\in\mathbb{N}$ for $0\le i,j\le m$ and $\mbox{vol}_{m}(\mathcal{S})\neq 0$ can be
  transformed via an isometry into the coordinates
  \begin{eqnarray*}
    v_0&=&(0,0,\dots,0),\\
    v_1&=&(q_{1,1}\sqrt{k_1},0,0\dots,0),\\
    v_2&=&(q_{2,1}\sqrt{k_1},q_{2,2}\sqrt{k_2},0,\dots,0),\\
    \vdots\\
    v_{m}&=&(q_{m,1}\sqrt{k_1},q_{m,2}\sqrt{k_2},\dots,q_{m,m}\sqrt{k_{m}}),
  \end{eqnarray*}
  where $k_i$ is the squarefree part of 
  $\frac{\mbox{vol}_i(v_0',v_1',\dots,v_i')^2}{\mbox{vol}_{i-1}(v_0',v_1',\dots,v_{i-1}')^2}$,
  $q_{i,j}\in\mathbb{Q}$, and $q_{j,j},k_j\neq 0$.
\end{Lemma}

We remark that we always have $k_1=1$. The connection between the $k_i$ and the characteristic $\mbox{char}(\mathcal{P})=k$ is given by
$$
  \mbox{char}(\mathcal{P})=\mbox{char}(\mathcal{S})=k=\text{square free part of }\prod_{i=1}^m k_i.
$$
Thus plane integral point sets $\mathcal{P}$ with rational coordinates are exactly those with characteristic $\mbox{char}(\mathcal{P})=1$. Due to Conjecture \ref{conj_embedding} plane integral point sets over $\mathbb{Z}^2$ correspond to plane integral point sets with characteristic $1$. So in principle there is no need to care about the coordinates -- this can still be done afterwards.

\medskip

There is one further transformation that maps integral point sets over $\mathbb{Z}^2$ onto integral point sets over $\mathbb{Z}^2$: scaling by an integral factor $\lambda$. One handicap of this mapping is that the inverse mapping may lead to non-integral point sets. Another shortcoming is that maximal integral point sets may be mapped onto non-maximal integral point sets. An example is given by the maximal integral point set $\mathcal{P}=\left\{
\begin{pmatrix}0\\0\\\end{pmatrix},
\begin{pmatrix}3\\0\\\end{pmatrix},
\begin{pmatrix}0\\4\\\end{pmatrix},
\begin{pmatrix}3\\4\\\end{pmatrix}
\right\}$.
If we scale it by a factor of $2$ we obtain $2\cdot\mathcal{P}=\left\{
\begin{pmatrix}0\\0\\\end{pmatrix},
\begin{pmatrix}6\\0\\\end{pmatrix},
\begin{pmatrix}0\\8\\\end{pmatrix},
\begin{pmatrix}6\\8\\\end{pmatrix}
\right\}$ an integral point set over $\mathbb{Z}^2$ which can be extended by the point $\begin{pmatrix}3\\4\end{pmatrix}$. In contrast to this example the integral point set $3\cdot\mathcal{P}=\left\{
\begin{pmatrix}0\\0\\\end{pmatrix},
\begin{pmatrix}9\\0\\\end{pmatrix},
\begin{pmatrix}0\\12\\\end{pmatrix},
\begin{pmatrix}9\\12\\\end{pmatrix}
\right\}$
is maximal. One might conjecture that for every maximal integral point set $\mathcal{M}$ there exists an integer $\lambda>1$ such that $\lambda\cdot\mathcal{M}$ is also maximal.

\section{Maximal integral point sets with given cardinality and minimum diameter}
\label{sec_maximal_integral_point_sets}

\noindent
From the combinatorial point of view a natural question is to ask for the minimum possible diameter $d_M(k,m)$ of a maximal integral point set $\mathcal{M}\subset\mathbb{Z}^m$ of cardinality $k$. If such a point set does not exist we set $d_M(k,m)=\infty$. Utilizing the exhaustive algorithm described in Section \ref{sec_exhaustive_generation} we have obtained the results given in Table \ref{table_min_diameter}.

\begin{table}[htp]
  \begin{center}
    \begin{tabular}{|r|r|r|}
      \hline
      $\mathbf{k}$ & $\mathbf{d_M(k,2)}$ & \textbf{corresponding point set} \\
      \hline
       4 &        $5$ & $\{(0,0),(3,4),(0,4),(3,0)\}$\\
      \hline
       5 &        $8$ & $\{(0,0),(3,4),(0,4),(0,8),(-3,4)\}$\\
      \hline
       6 &       $25$ & $\{(0,0),(12,16),(12,9),(-12,9),(-12,16),(0,25)\}$\\
      \hline
       7 &       $30$ & $\{(0,0),(6,8),(0,8),(0,16),(-6,8),(-15,8),(15,8)\}$\\
      \hline
       8 &       $65$ & $\{(0,0),(15,36),(0,16),(15,-20),(48,-20),(48,36),(63,0),(63,16)\}$\\
      \hline
       9 &       $96$ & $\{(0,0),(15,20),(0,20),(0,40),(0,56),(0,-16),(-15,20),(-48,20),(48,20)\}$\\
      \hline
         &            & $\{(0,0),(22,120),(0,120),(-27,120),(160,120),(182,0),(182,120),$\\
      10 &  $\le 600$ & $(-209,120),(209,120),(391,120)\}$\\
      \hline
         &            & $\{(0,0),(5,12),(0,12),(0,24),(-5,12),(-9,12),(9,12),(-16,12),(16,12),$\\
      11 &       $70$ & $(-35,12),(35,12)\}$\\
      \hline
         &            & $\{(0,0),(35,120),(35,84),(-64,-48),(0,204),(-189,-48),(-64,252),$\\
      12 &  $\le 325$ & $(-253,0),(-189,252),(-288,84),(-288,120),(-253,204)\}$\\
      \hline
         &            & $\{(0,0),(48,64),(0,64),(0,128),(-48,64),(-120,64),(120,64),(-252,64),$\\
      13 & $\le 2046$ & $(252,64),(-510,64),(510,64),(-1023,64),(1023,64)\}$\\
      \hline
    \end{tabular}
    \caption{Minimum possible diameters of maximal plane integral point sets with given cardinality.}
    \label{table_min_diameter}
  \end{center}
\end{table}

Clearly we have $d_M(1,2)=d_M(2,2)=\infty$ since a line $l$ through two different points $P_1$ and $P_2$ with integral coordinates and integral distance $\overline{P_1P_2}$ contains an infinite integral point set $\mathcal{P}=\{P_1+\lambda\cdot(P_2-P_1)\,:\,\lambda\in\mathbb{Z}\}$ as a subset. So the next value to determine is $d_M(3,2)$. Whether $d_M(3,2)$ is finite had been an open question of \cite{gauss_integers}, which was answered in \cite{paper_axel} by determining $d_M(3,2)=2066$, -- a diameter out of reach for our general exhaustive algorithm described in Section \ref{sec_exhaustive_generation}. But it can be easily adapted for this purpose. We alter Algorithm \ref{algo_triangles} by omitting  right-angled triangles, since these obviously are not maximal. Then we skip Algorithm \ref{algo_embedding} and directly run the version of Algorithm \ref{algo_enlargement} where we search for rational instead of integral solutions $(x_P,y_P)$ of equation system (\ref{eq_fourth_point}). If we have found the first solution $(x_P,y_P)$ for a given triangle $\Delta$ we can immediately stop our investigations on $\Delta$ since it cannot be a maximal integral triangle. Using these reductions and skipping the time consuming clique search we were able to exhaustively search for (strongly) maximal integral triangles over $\mathbb{Z}^2$ with diameter at most $15000$ \cite{paper_axel,herontriangles}. There are exactly $126$ such examples. Here we list the first, with respect to their diameter, ten examples, where we give the edge lengths and the coordinates in minimal list representation, which is unique in these cases:
\begin{eqnarray*}
 \{2066,1803,505\}  & & \left[(0,0)^T,\,(-336,-377)^T,\,(384,-2030)^T\,\right]    \\
 \{2549,2307,1492\} & & \left[(0,0)^T,\,(-700,-2451)^T,\,(1100,-1008)^T\,\right]  \\
 \{3796,2787,2165\} & & \left[(0,0)^T,\,(-387,-2760)^T,\,(1680,-3404)^T\,\right]  \\
 \{4083,2425,1706\} & & \left[(0,0)^T,\,(-410,-1656)^T,\,(1273,2064)^T\,\right]   \\
 \{4426,2807,1745\} & & \left[(0,0)^T,\,(-280,-2793)^T,\,(376,-4410)^T\,\right]   \\
 \{4801,2593,2210\} & & \left[(0,0)^T,\,(-1488,-1634)^T,\,(1632,2015)^T\,\right]  \\
 \{4920,4177,985\}  & & \left[(0,0)^T,\,(-473,-864)^T,\,(4015,1152)^T\,\right]    \\
 \{5044,4443,2045\} & & \left[(0,0)^T,\,(-1204,-1653)^T,\,(2156,-4560)^T\,\right] \\
 \{5045,4803,244\}  & & \left[(0,0)^T,\,(-44,-240)^T,\,(240,4797)^T\,\right]      \\
 \{5186,5163,745\}  & & \left[(0,0)^T,\,(-407,-624)^T,\,(4030,-3264)^T\,\right]   \\
\end{eqnarray*}

\section{Constructions for maximal integral point sets over $\mathbb{Z}^2$}
\label{sec_constructions}

\noindent
In this section we want to describe constructions for maximal integral point sets $\mathcal{M}$ of a given cardinality or a given shape. In most cases our constructions do not lead to integral point sets which are maximal in every case, but which yield candidates, which are very likely to be maximal (from an empiric point of view). W.l.o.g. we can assume that the origin $(0,0)^T$ is always contained in $\mathcal{M}$. Every further point $(a,b)^T$ meets $a^2+b^2=c^2$. In this case we call $(a,b)$ a Pythagorean pair or $(a,b,c)$ a Pythagorean triple. If additionally $\gcd(a,b)=\gcd(a,b,c)=1$ we speak of primitive pairs or triples. Given only one Pythagorean pair $(a,b)$ we can perform the following two constructions for integral point sets over $\mathbb{Z}^2$:

\begin{Construction}
  \label{construction_1}
  If $(a,b)$ is a Pythagorean pair, then $\mathcal{P}_1(a,b):=\left\{
  \begin{pmatrix}0\\0\end{pmatrix},
  \begin{pmatrix}a\\0\end{pmatrix},
  \begin{pmatrix}0\\b\end{pmatrix},
  \begin{pmatrix}a\\b\end{pmatrix}
  \right\}$ is an integral point set of cardinality $4$.
\end{Construction}

\begin{Construction}
  \label{construction_2}
  If $(a,b)$ is a Pythagorean pair, then $\mathcal{P}_2(a,b):=\left\{
  \begin{pmatrix}0\\0\end{pmatrix},
  \begin{pmatrix}a\\0\end{pmatrix},
  \begin{pmatrix}-a\\0\end{pmatrix},
  \begin{pmatrix}0\\b\end{pmatrix},
  \begin{pmatrix}0\\-b\end{pmatrix}
  \right\}$ is an integral point set of cardinality $5$.
\end{Construction}

\begin{figure}[htp]
  \begin{center}
    \setlength{\unitlength}{0.18cm}
    \begin{picture}(70,24)
      \put(35,12){\circle*{1}}
      \put(40,12){\circle*{1}}%
      \put(35,24){\circle*{1}}
      \put(35,0){\circle*{1}}
      \put(30,12){\circle*{1}}
      \put(44,12){\circle*{1}}
      \put(26,12){\circle*{1}}
      \put(51,12){\circle*{1}}
      \put(19,12){\circle*{1}}
      \put(70,12){\circle*{1}}
      \put(0,12){\circle*{1}}
      \put(0,12){\line(1,0){70}}
      \put(35,24){\line(0,-1){24}}
      \qbezier(35,24)(37.5,18)(40,12)
      \qbezier(35,24)(32.5,18)(30,12)
      \qbezier(35,24)(39.5,18)(44,12)
      \qbezier(35,24)(30.5,18)(26,12)
      \qbezier(35,24)(43,18)(51,12)
      \qbezier(35,24)(27,18)(19,12)
      \qbezier(35,24)(52.5,18)(70,12)
      \qbezier(35,24)(17.5,18)(0,12)
      \qbezier(35,0)(37.5,6)(40,12)
      \qbezier(35,0)(32.5,6)(30,12)
      \qbezier(35,0)(39.5,6)(44,12)
      \qbezier(35,0)(30.5,6)(26,12)
      \qbezier(35,0)(43,6)(51,12)
      \qbezier(35,0)(27,6)(19,12)
      \qbezier(35,0)(52.5,6)(70,12)
      \qbezier(35,0)(17.5,6)(0,12)
    \end{picture}
  \end{center}
  \caption{Extension of $\mathcal{P}_2(5,12)$ to a crab of cardinality $11$.}
  \label{fig_crab_11}
\end{figure}
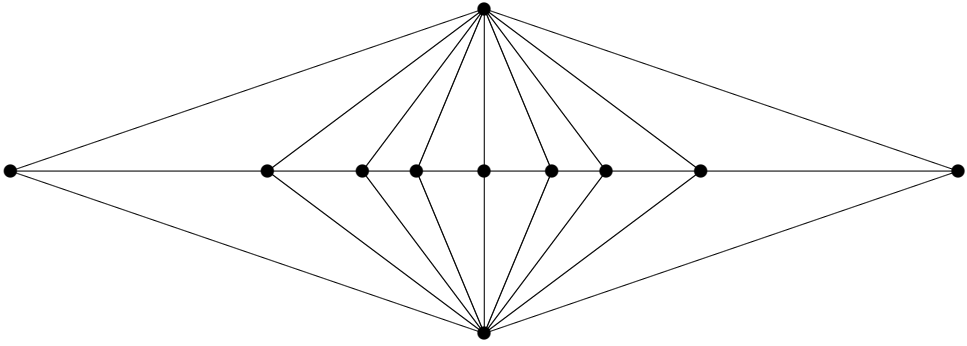

We call Construction \ref{construction_1} the \textit{rectangle construction} of $(a,b)$ and Construction \ref{construction_2} the \textit{rhombus construction} of $(a,b)$. If we choose $(a,b)$ with $2|a$, $2|b$ then clearly $\mathcal{P}_1(a,b)$ cannot be maximal. On the other side $\mathcal{P}_1(9,12)$ is a maximal integral point set although $\gcd(9,12)=3$. Empirically, we have observed that for primitive pairs $(a,b)$ the point set $\mathcal{P}_1(a,b)$ is maximal in many, but not all cases, see e.g. the non maximal integral point set $\mathcal{P}_1(7,24)$, which can be extended to the maximal integral point set $\left\{
\begin{pmatrix}0\\0\\\end{pmatrix},
\begin{pmatrix}7\\0\end{pmatrix},
\begin{pmatrix}0\\24\end{pmatrix},
\begin{pmatrix}7\\24\end{pmatrix},
\begin{pmatrix}-9\\12\end{pmatrix},
\begin{pmatrix}16\\12\end{pmatrix}
\right\}$. For $(a,b)=(3,4)$ both constructions $\mathcal{P}_1(a,b)$ and $\mathcal{P}_2(a,b)$ yield maximal integral point sets. Empirically Construction \ref{construction_2} is a bit weaker, since it often happens that $\mathcal{P}_1(a,b)$ is maximal but $\mathcal{P}_2(a,b)$ is not, as for example for $(a,b)=(5,12)$. For the other direction we have no example.
We would like to mention that $\mathcal{P}_2(5,12)$ can be extended to the very interesting maximal integral point set $\mathcal{M}=\left\{
\begin{pmatrix}0\\0\end{pmatrix},
\begin{pmatrix}5\\0\end{pmatrix},
\begin{pmatrix}0\\12\end{pmatrix},
\begin{pmatrix}0\\-12\end{pmatrix},
\begin{pmatrix}-5\\0\end{pmatrix},
\begin{pmatrix}9\\0\end{pmatrix},
\begin{pmatrix}-9\\0\end{pmatrix},
\begin{pmatrix}16\\0\end{pmatrix},\right.$
$\left.
\begin{pmatrix}-16\\0\end{pmatrix},
\begin{pmatrix}35\\0\end{pmatrix},
\begin{pmatrix}-35\\0\end{pmatrix}
\right\}$, which has an intriguing geometrical structure, see Figure \ref{fig_crab_11}.

\begin{Definition}
  For positive integers $a,b_1,\dots,b_k$ we call the point set
  $$
    \mbox{crab}(a,b_1,\dots,b_k):=\left\{
    \begin{pmatrix}0\\0\end{pmatrix},
    \begin{pmatrix}0\\\pm a\end{pmatrix},
    \begin{pmatrix}\pm b_1\\0\end{pmatrix},
    \dots,
    \begin{pmatrix}\pm b_k\\0\end{pmatrix}
    \right\}
  $$
  a crab of order $k$.
\end{Definition}
We remark that the cardinality of $\mbox{crab}(a,b_1,\dots,b_k)$ is given by $2k+3$ and that the point set is symmetric w.r.t. the two coordinate axes. This point set is indeed integral if the pairs $(a,b_1),\dots,(a,b_k)$ are Pythagorean pairs. So it is very easy to construct crabs, either directly or by extending $\mathcal{P}_2(a,b)$, see Subsection \ref{subsec_crabs}. Empirically the extension points of $\mathcal{P}_2(a,b)$ very often lie on one of the two axis. An example that this must not be the case in general is given by the primitive pair $(1480,969)$, where $\mathcal{P}_2(1480,969)$ can be extended to $\left\{
\begin{pmatrix}0\\0\end{pmatrix},
\begin{pmatrix}1480\\0\end{pmatrix},
\begin{pmatrix}-1480\\0\end{pmatrix},
\begin{pmatrix}0\\969\end{pmatrix},
\begin{pmatrix}0\\-969\end{pmatrix},
\begin{pmatrix}1040\\462\end{pmatrix},
\begin{pmatrix}1040\\-462\end{pmatrix},
\begin{pmatrix}-1040\\462\end{pmatrix},\right.$
$\left.\begin{pmatrix}-1040\\-462\end{pmatrix}
\right\}$.

\subsection{Construction of crabs}
\label{subsec_crabs}

\begin{figure}[ht]
  \begin{center}
    \setlength{\unitlength}{0.7cm}
    \begin{picture}(9,4.5)
      \put(0,0){\line(1,0){9}}
      \put(0,0){\line(5,4){5}}
      \put(1,0){\line(1,1){4}}
      \put(2,0){\line(3,4){3}}
      \put(4,0){\line(1,4){1}}
      \put(6,0){\line(-1,4){1}}
      \put(7,0){\line(-1,2){2}}
      \put(9,0){\line(-1,1){4}}
      \put(0,0){\circle*{0.3}}
      \put(1,0){\circle*{0.3}}
      \put(2,0){\circle*{0.3}}
      \put(4,0){\circle*{0.3}}
      \put(6,0){\circle*{0.3}}
      \put(7,0){\circle*{0.3}}
      \put(9,0){\circle*{0.3}}
      \put(5,4){\circle*{0.3}}
      \put(5,0){\line(0,1){4}}
      \put(0.1,-0.5){$a_3$}
      \put(1.1,-0.5){$a_2$}
      \put(2.8,-0.5){$a_1$}
      \put(4.4,0.2){$q$}
      \put(5.2,0.2){$q'$}
      \put(5.1,1.4){$h$}
      \put(4.7,-0.5){$a_0$}
      \put(6.2,-0.5){$a_1'$}
      \put(7.6,-0.5){$a_2'$}
      \put(1.6,2){$b_3$}
      \put(0.8,0.4){$b_2$}
      \put(1.7,0.4){$b_1$}
      \put(3.5,0.4){$b_0$}
      \put(6.05,0.4){$b_0'$}
      \put(7.0,0.4){$b_1'$}
      \put(7,2){$b_2'$}
    \end{picture}\\[2mm]
    \caption{Plane integral point set $\mathcal{P}$ with $n-1$ points on a line.}
    \label{fig_on_a_line}
  \end{center}
\end{figure}
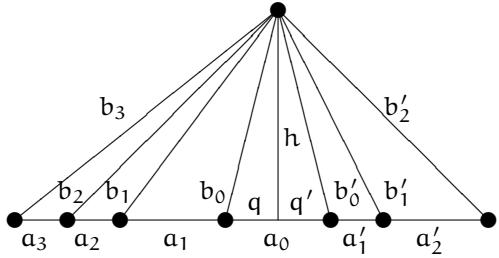

Since many maximal integral point sets over $\mathbb{Z}^2$ are crabs we are interested in a method to construct them directly. From the general theory of integral point sets we know that integral point sets $\mathcal{P}$ over $\mathbb{R}^2$ with minimum diameter consist of point sets with $n-1$ collinear points, see Figure \ref{fig_on_a_line}, for $9\le n\le 122$ points, see \cite{phd_kurz,paper_alfred}. For these point sets there is an interesting connection between the points of the point set $\mathcal{P}$ and divisors of a certain number $D$, see \cite{phd_kurz,paper_alfred}.

\begin{Definition}
  The decomposition number $D$ of an integral triangle with side lengths $a$, $b$, and $c$ is given by
  $$
    D=\frac{(a+b+c)(a+b-c)(a-b+c)(-a+b+c)}{\gcd(b^2-c^2+a^2,2a)^2}.
  $$
\end{Definition}

\begin{Lemma}{\bf(Decomposition lemma)\\}
  \label{lemma_decomposition}
  \noindent
  The distances of a plane integral point set $\mathcal{P}$ consisting of $n$ points where a subset of $n-1$ points is 
  collinear correspond to decompositions of the decomposition number $D$ of the largest triangle of $\mathcal{P}$ 
  into two factors.
\end{Lemma}
\begin{proof}
  We use the notation of Figure \ref{fig_on_a_line} and set 
  $$c_i=q+\sum_{j=1}^ia_j\quad\mbox{for}\quad 0\le i\le s,\, c_i'=q'+\sum_{j=1}^ia_j'\quad\mbox{for}\quad 0\le i\le t.$$
  Pythagoras' Theorem yields $c_{i+1}^2+h^2=b_{i+1}^2$ and $c_i^2+h^2=b_i^2$ for $0\le i< s$. We subtract these equations from 
  each other and get $$a_{i+1}^2+2a_{i+1}\sum_{j=1}^ia_j+2a_{i+1}q=b_{i+1}^2-b_i^2\,.$$ Because the $a_i$ and the $b_i$ are	
  positive integers we have $2a_{i+1}q\in\mathbb{N}$ for $0\le i<s$ and therefore $2\gcd(a_1,a_2,\dots,a_s)q\in\mathbb{N}$. From 
  $q+q'=a_0\,\,\in\mathbb{N}$ we conclude $2\gcd(a_1,a_2,\dots,a_s)q'\in\mathbb{N}$. With an analogous conclusion for the $c_i'$ 
  and $g=2\gcd(a_1,\dots,a_s,a_1',\dots,a_t')$ we get $gq\in\mathbb{N}\quad\mbox{and}\quad gq'\in\mathbb{N}$. A last use of 
  Pythagoras' Theorem yields for $1\le i\le s$ and for $1\le j\le t$ the factorization of $g^2h^2$ into a product of two positive
  integers,
  $$g^2h^2=(gb_i+gc_i)(gb_i-gc_i)=(gb_j'+gc_j')(gb_j'-gc_j').$$
  So we can obtain the possible values for $c_i$ and $c_i'$ by decomposing $g^2h^2$ into two factors.

  If we are given the three side lengths $a$, $b$, and $c$ of an integral triangle and want to determine the points on the side
  of length $a$ so that the resulting point set is integral then we can associate $b$ with $b_s$, $c$ with $b_t'$, and $a$ with
  $\sum_{i=1}^sa_i+a_0+\sum_{i=1}^ta_i'$. With this we have $$c_s=q+\sum\limits_{j=1}^sa_j=\frac{b^2-c^2+a^2}{2a}\,.$$
  Because $g$ can also be defined as the smallest integer with $gc_s\in\mathbb{N}$ we receive
  $$g=\frac{2a}{\gcd(b^2-c^2+a^2,2a)}\,.$$
  Due to the Heron formula 
  $16A_\Delta^2=(a+b+c)(a+b-c)(a-b+c)(-a+b+c)$
  and the formula for the area of a triangle
  $2A_\Delta=ah$
  we finally get
  \begin{eqnarray*}
    g^2h^2=\frac{g^2(a+b+c)(a+b-c)(a-b+c)(-a+b+c)}{4a^2}=\\
    =\frac{(a+b+c)(a+b-c)(a-b+c)(-a+b+c)}{\gcd(b^2-c^2+a^2,2a)^2}=D\,.
  \end{eqnarray*}
\end{proof}

If we choose $g=1$ and $h\in\mathbb{N}$ we can directly apply Lemma \ref{lemma_decomposition} to construct crabs. Let us look at an example. We choose $g=1$ and $h=2\cdot3\cdot 5=30$. The divisors of $D=g^2h^2=900$ are given by $\{1,2,3,4,5,6,9,10,12,15,18,20,25,30,36,45,50,60,75,90,100,150,180,225,300,450,900\}$. If we have $D=f_1\cdot f_2$, then $b=\frac{f_1+f_2}{2}$ and $c=\frac{f_1-f_2}{2}$. Thus we must have $f_1>f_2$ and $f_1\equiv f_2\pmod 2$ to determine the values $b_i$ of a corresponding crab. Here we have $b_1=\frac{50-18}{2}=16$, $b_2=\frac{90-10}{2}=40$, $b_3=\frac{150-6}{2}=72$,  $b_4=\frac{450-2}{2}=224$, and $a=h=30$. This yields the integral point set $\mbox{crab}(30,16,40,72,224)$ of cardinality $2\cdot 4+3=11$ and diameter $\max\left\{2b_i,2a,\sqrt{b_i^2+a^2}\right\}=2\cdot\max\{b_i,a\}=448$. Given the prime factorization $h=\prod_{i=1}^r p_i^{\alpha_i}$ it is not difficult to determine the $k$-value of the resulting crab. Let us fix $p_1=2$ and set $\tilde{\alpha}_1=\max\left(\alpha_1-1,0\right)$. With this we can state
\begin{equation}
  \label{eq_crab_order}
  k=\frac{(2\tilde{\alpha}_1+1)\cdot\prod_{i=2}^r (2\alpha_i+1)\,-1}{2}.
\end{equation}
Using $h=p^k$, where $p$ is an arbitrary odd prime, we are able to produce a crab of order $k$ for each $k\ge 1$. Thus we have constructions for integral point sets of cardinality $2k+3$ for each $k\in\mathbb{N}$. To obtain small point sets with many points we should clearly choose integers with many divisors for $h$ instead. As for all of our constructions the maximality of the resulting integral point set is not guaranteed, but very likely.

\begin{Construction}
  \label{line_construction}
  For a given integer $h$ there exists an integral point set $\mbox{decompose}(h)$ which is a crab of order $k$, where $k$ is
  given by Equation (\ref{eq_crab_order}).
\end{Construction}

If $h>4$ then the diameter of $\mbox{decompose}(h)$ is given by $h^2-1$ if $h$ is odd and given by $\frac{h^2}{2}-2$ if $h$ is even.

\begin{Conjecture}
  For each integer $h$ the plane integral point set $\mathcal{P}=\mbox{decompose}(h)$ is maximal if $|\mathcal{P}|\ge 7$.
\end{Conjecture}

Also, the recognition of a crab is a very easy task. Given an integral point set $\mathcal{P}$ over $\mathbb{Z}^2$ one can easily check whether a subset $\mathcal{L}\subset\mathcal{P}$ of $n-2$ points is collinear by using:
\begin{Lemma}
  \label{lemma_collinear}
  Three points $(x_1,y_1)$, $(x_2,y_2)$, and $(x_3,y_3)$ in $\mathbb{R}^2$ are collinear if and only if we have
  $$
    \left|
    \begin{array}{rrr}
      x_1 & y_1 & 1 \\
      x_2 & y_2 & 1 \\
      x_3 & y_3 & 1 \\
    \end{array}
    \right|=0.
  $$
\end{Lemma}
Additionally the lines through $\mathcal{L}$ and $\mathcal{P}\backslash\mathcal{L}$ are perpendicular. If the point set is symmetric to these two lines then $\mathcal{P}$ is a crab.

Crabs are very dominating examples of maximal integral point sets over $\mathbb{Z}^2$. For the study of maximal integral point sets over $\mathbb{Z}^3$ one might try to generalize the construction of a crab. Let us remark in this context, that the existence of an integral point set with coordinates
$$
  \left\{
  \begin{pmatrix}0\\0\\0\\\end{pmatrix},
  \begin{pmatrix}x\\0\\0\\\end{pmatrix},
  \begin{pmatrix}0\\y\\0\\\end{pmatrix},
  \begin{pmatrix}0\\0\\z\\\end{pmatrix}
  \right\},
$$
where $x,y,z\in\mathbb{Z}$ is equivalent to a famous open problem, the existence of a perfect box, see \cite[D18]{UPIN}.

So far we have only used $g=1$ in Lemma \ref{lemma_decomposition}. Now we want to have a look at the case $g>1$. So given integers $g\cdot h$ and $g$ we can apply Lemma \ref{lemma_decomposition}. For two factors $f_1>f_2$ with $f_1\cdot f_2=g^2h^2$ we have
$$
  gb_i=\frac{f_1+f_2}{2}\text{ and }gc_i=\frac{f_1-f_2}{2}.
$$
The values $gb_i$ and $gc_i$ are integers if and only if we have $f_1\equiv f_2\pmod 2$. Since not only the $gb_i$'s but also the $b_i$'s must be integers we have to require $f_1+f_2\equiv 0\pmod g$. Let us have an example. We choose $gh=672=2^5\cdot 3\cdot 7$ and $g=5$. Now we look at the divisors of $g^2h^2=451584=2^{10}\cdot 3^3\cdot 7^2$ and determine the suitable pairs $(f_1,f_2)$ fulfilling
$$
  f_1\cdot f_2=g^2h^2,\quad
  f_1>f_2,\quad
  f_1\equiv f_2\pmod 2,\text{ and }
  f_1+f_2\equiv 0\pmod 5,
$$
\begin{eqnarray*}
  \!\!\!\!\!\!\!&&\Big\{
    (784,576),
    (896,504),
    (1176,384),
    (1344,336),
    (1536,294),
    (1764,256),
    (2016,224),\\
    \!\!\!\!\!\!\!&&(2304,196),
    (3136,144),
    (3584,126),
    (4704,96),
    (5376,84),
    (7056,64),
    (8064,56),
    (12544,36),\\
    \!\!\!\!\!\!\!&&(18816,24),
    (28224,16),
    (32256,14),
    (75264,6),
    (112896,4)
  \Big\}.
\end{eqnarray*}
The corresponding values $b_i$ are given by
$$
  \left\{136,140,156,168,183,202,224,250,328,371,480,546,712,812,1258,1884,2824,3227,7527,11290\right\}
$$
and the corresponding values $gc_i$ are given by
\begin{eqnarray*}
  &&\mathcal{C}=\Big\{104,196,396,504,621,754,896,1054,1496,1729,2304,\\
  &&2646,3496,4004,6254,9396,14104,16121,37629,56446\Big\}.
\end{eqnarray*}
Clearly the $c_i$ cannot be integers unless $h$ is an integer. So let us consider the points on the left of the base point $F$ of the height $h$. They correspond to values $gc_i$ which all fulfill $gc_i\equiv m\pmod g$, for a fixed $m\in\{1,\dots,g-1\}$. The points on the right hand side of $F$ correspond to the values $gc_i$ fulfilling $gc_i\equiv -m\pmod g$. So let us choose $m=1$. Since all elements of our candidate set $\mathcal{C}$ are congruent to $\pm 1$ modulo $5$ we obtain an integral point set of cardinality $|\mathcal{C}|+1=21$:
\begin{eqnarray*}
  \mathcal{P} &=& \Big\{
  \begin{pmatrix}0\\\frac{672}{5}\end{pmatrix},
  \begin{pmatrix}\frac{-196}{5}\\0\end{pmatrix},
  \begin{pmatrix}\frac{-396}{5}\\0\end{pmatrix},
  \begin{pmatrix}\frac{-621}{5}\\0\end{pmatrix},
  \begin{pmatrix}\frac{-896}{5}\\0\end{pmatrix},
  \begin{pmatrix}\frac{-1496}{5}\\0\end{pmatrix},
  \begin{pmatrix}\frac{-2646}{5}\\0\end{pmatrix},
  \begin{pmatrix}\frac{-3496}{5}\\0\end{pmatrix},\\
  &&\begin{pmatrix}\frac{-9396}{5}\\0\end{pmatrix},
  \begin{pmatrix}\frac{-16121}{5}\\0\end{pmatrix},
  \begin{pmatrix}\frac{-56446}{5}\\0\end{pmatrix},
  \begin{pmatrix}\frac{104}{5}\\0\end{pmatrix},
  \begin{pmatrix}\frac{504}{5}\\0\end{pmatrix},
  \begin{pmatrix}\frac{754}{5}\\0\end{pmatrix},
  \begin{pmatrix}\frac{1054}{5}\\0\end{pmatrix},
  \begin{pmatrix}\frac{1729}{5}\\0\end{pmatrix},\\
  &&\begin{pmatrix}\frac{2304}{5}\\0\end{pmatrix},
  \begin{pmatrix}\frac{4004}{5}\\0\end{pmatrix},
  \begin{pmatrix}\frac{6254}{5}\\0\end{pmatrix},
  \begin{pmatrix}\frac{14104}{5}\\0\end{pmatrix},
  \begin{pmatrix}\frac{37629}{5}\\0\end{pmatrix}.
  \Big\}
\end{eqnarray*}

After a suitable transformation and applying Algorithm \ref{algo_normalization} we obtain the minimum coordinate representation
\begin{eqnarray*}
  \Big[
  \begin{pmatrix}0\\0\end{pmatrix},
  \begin{pmatrix}0\\-168\end{pmatrix},
  \begin{pmatrix}-40\\30\end{pmatrix},
  \begin{pmatrix}64\\-48\end{pmatrix},
  \begin{pmatrix}-88\\66\end{pmatrix},
  \begin{pmatrix}112\\-84\end{pmatrix},
  \begin{pmatrix}144\\-108\end{pmatrix},
  \begin{pmatrix}180\\-135\end{pmatrix},\\
  \begin{pmatrix}-196\\147\end{pmatrix},
  \begin{pmatrix}224\\-168\end{pmatrix},
  \begin{pmatrix}-288\\216\end{pmatrix},
  \begin{pmatrix}320\\-240\end{pmatrix},
  \begin{pmatrix}504\\-378\end{pmatrix},
  \begin{pmatrix}-560\\420\end{pmatrix},
  \begin{pmatrix}640\\-480\end{pmatrix},\\
  \begin{pmatrix}-920\\690\end{pmatrix},
  \begin{pmatrix}1584\\-1188\end{pmatrix},
  \begin{pmatrix}-2176\\1632\end{pmatrix},
  \begin{pmatrix}2660\\-1995\end{pmatrix},
  \begin{pmatrix}-5940\\4455\end{pmatrix},
  \begin{pmatrix}9112\\-6834\end{pmatrix}
  \Big]
\end{eqnarray*}

\begin{figure}[htp]
  \begin{center}
    \setlength{\unitlength}{0.0008cm}
    \begin{picture}(16000,11300)
      \put(5940,6834){\circle*{120}}
      \put(5940,6666){\circle*{120}}
      \put(5900,6864){\circle*{120}}
      \put(6004,6786){\circle*{120}}
      \put(5852,6900){\circle*{120}}
      \put(6052,6750){\circle*{120}}
      \put(6084,6726){\circle*{120}}
      \put(6120,6699){\circle*{120}}
      \put(5744,6981){\circle*{120}}
      \put(6164,6666){\circle*{120}}
      \put(5652,7050){\circle*{120}}
      \put(6260,6594){\circle*{120}}
      \put(6444,6456){\circle*{120}}
      \put(5380,7254){\circle*{120}}
      \put(6580,6354){\circle*{120}}
      \put(5020,7524){\circle*{120}}
      \put(7524,5646){\circle*{120}}
      \put(3764,8466){\circle*{120}}
      \put(8600,4839){\circle*{120}}
      \put(0,11289){\circle*{120}}
      \put(15052,0){\circle*{120}}
      \qbezier(15052,0)(7526,5644.5)(0,11289)
    \end{picture}
    \caption{A semi-crab of cardinality $21$ and diameter $18815$.}
    \label{fig_semi_crab}
  \end{center}
\end{figure}
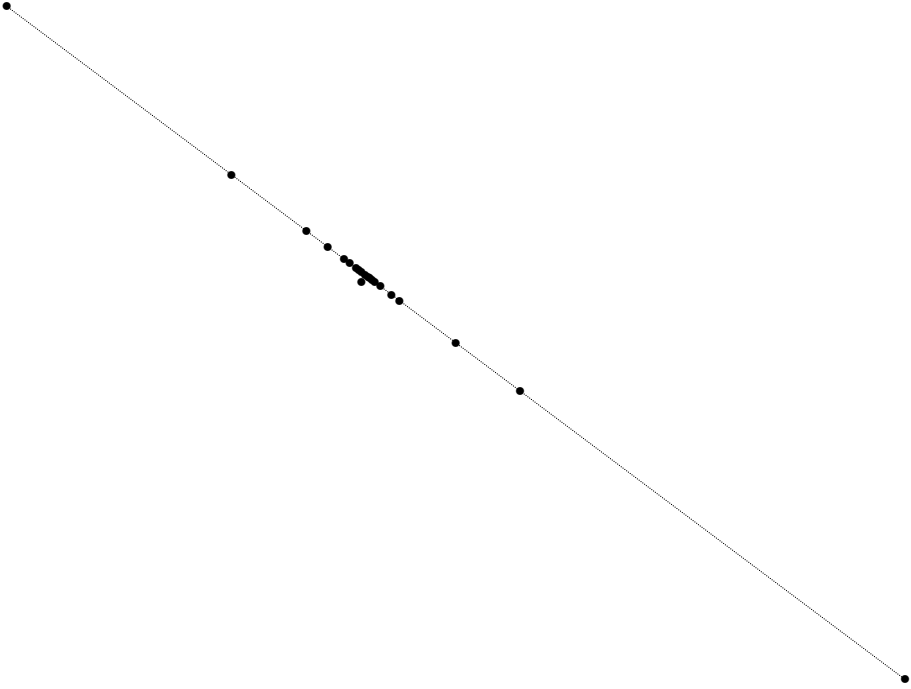

We call point sets arising from Lemma \ref{lemma_decomposition}, where $g>1$ and $h\notin\mathbb{N}$ semi-crabs, see Figure \ref{fig_semi_crab} for a drawing of our example.

\begin{Definition}
  For positive integers $g,\tilde{h}$ and non-zero integers $\tilde{b}_1,\dots,\tilde{b}_k$ we call the point set
  $$
    \mbox{semi-crab}\left(g,\tilde{h},\tilde{b}_1,\dots,\tilde{b}_k\right):=\left\{
    \begin{pmatrix}0\\\frac{h}{g}\end{pmatrix},
    \begin{pmatrix}\frac{\tilde{b}_1}{g}\\0\end{pmatrix},
    \dots,
    \begin{pmatrix}\frac{\tilde{b}_k}{g}\\0\end{pmatrix}
    \right\}
  $$
  a semi-crab of order $k$.
\end{Definition}

\begin{Construction}
  \label{semi_crab_construction}
  For given positive integers $g$ and $gh$, where $h\notin\mathbb{N}$, there exists an integral point set
  $\mbox{decompose}(g,gh)$ which is isomorphic to a semi-crab.
\end{Construction}

\begin{Conjecture}
  For each pair of integers $gh$, $g$ the plane integral point set $\mathcal{P}=\mbox{decompose}(gh,g)$ is maximal if
  $|\mathcal{P}|\ge 7$.
\end{Conjecture}

Since Construction \ref{line_construction} can only produce integral point sets of odd cardinality Construction \ref{semi_crab_construction} is a convenient completion. It is not difficult to implement an algorithm that exhaustively generates crabs and semi-crabs up to a given diameter by utilizing Lemma \ref{lemma_decomposition}.

Let us have a look at the possible values for $g>1$. If we would choose $g=2$ then due to $2\nmid g^2h^2$ all divisors of $g^2h^2$ would be odd and we would have $m\equiv 1\pmod 2$, which is not possible. Thus $2\nmid g$. For $g=3$ the only possibility for $f_1+f_2\equiv 0\pmod 3$ is $f_1\equiv 1\pmod 3$ and $f_2\equiv 2\pmod 3$, which is not possible since $1\cdot 2\equiv 2\pmod 3$ is not a square in $\mathbb{Z}_3$. Thus $g=5$ is the first valid possibility. More generally we can state that if $g$ is a prime then we have $g\equiv 1\pmod 4$, since $-1$ has to be a square in $\mathbb{Z}_g$.

\subsection{Construction of integral point sets on circles}
\label{subsec_circles}

In Addition to the construction of crabs there is another useful construction of integral point sets of $\mathbb{Z}^2$ with large cardinality, see \cite{minimum_diameter} for a similar construction over the ring $\mathbb{Z}\left[\frac{-1+\sqrt{-3}}{2}\right]$. Let $p_j\equiv 1\pmod 4$ be distinct primes over $\mathbb{N}$. We consider the ring $\mathbb{Z}[i]$, where every integer $p_j$ has a unique prime factorization $p_j=\omega_j\cdot\overline{\omega_j}$, where $\overline{c}$ denotes the complex conjugate of $c$. We may write $\omega_{j}=a_j+b_ji$, with integers $a_j,b_j$. With multiplicities $v_j\in\mathbb{N}$ we set
$$
  R=\prod_{j=1}^r p_j^{v_j}
$$
and for each of the $\tau(R)$ divisors of $R$,
$$
  \prod_{j=1}^r p_j^{u_j}=\prod_{j=1}^r\omega_j^{u_j}\overline{\omega}_j^{u_j},\quad 0\le u_j\le v_j
$$
we define $\eta_{2h}=\prod_{j=1}^r\omega_j^{v_j+u_j}\overline{\omega}_j^{v_j-u_j}$, $\eta_{2h-1}=i\cdot\eta_{2h}$ for $1\le h\le\tau(R)$. With this we define vertices $\xi_s$ for $1\le s\le 2\tau(R)$ by
$$
  \xi_{2h-k}=\frac{\eta_{2h-k}^2}{R},\quad 1\le h\le\tau(R),\, k\in\{0,1\}.
$$
We set $\eta_s=x_s+y_si$ with $x_s,y_s\in\mathbb{Z}$ for $1\le s\le 2\tau(R)$. We have
$$
  |\eta_s|^2=\eta_s\overline{\eta}_s=x_s^2+y_s^2=\prod_{j=1}^r\omega_j^{2v_j}\overline{\omega}_j^{2v_j}=\prod_{j=1}^r p_j^{2v_j}=R^2.
$$
This yields $x_s^2=R^2-y_s^2$, which we use to calculate
\begin{eqnarray*}
  R^2\cdot|\xi_s-\xi_t|^2 &=& |\eta_s^2-\eta_t^2|^2\\
                          &=& |x_s^2-y_s^2-x_t^2+y_t^2+i\cdot(2x_sy_s-2x_ty_t)|^2\\
                          &=& |(2y_t^2-2y_s^2)+i\cdot(2x_sy_s-2x_ty_t)|^2\\
                          &=& 4(y_t^2-y_s^2)^2+4(x_sy_s-x_ty_t)^2\\
                          &=& 4(y_t^2-y_s^2)(x_s^2-x_t^2)+4(x_sy_s-x_ty_t)^2\\
                          &=& 2^2(x_sy_t-x_ty_s)^2.
\end{eqnarray*}
Thus the distance between $\xi_s$ and $\xi_t$ is given by $|\xi_s-\xi_t|=\frac{1}{R}|x_sy_t-x_ty_s|$. Since
$$
  \eta_s\overline{\eta}_t=(x_s+y_si)(x_t-y_ti)=x_sx_t+y_sy_t+i(x_ty_s-x_sy_t)
$$
and
\begin{eqnarray*}
  \eta_s\overline{\eta}_t&=&i^{k_s}\overline{i}^{k_t}\prod_{j=1}^r \omega_j^{v_j+u_j}\overline{\omega}_j^{v_j-u_j}
  \prod_{j=1}^r \overline{\omega}_j^{v_j+w_j}\omega_j^{v_j-w_j}\\
  &=& i^{k_s-k_t}\prod_{j=1}^r\omega_j^{2v_j+u_j-w_j}\overline{\omega}_j^{2v_j-u_j+w_j}\\
  &=& R\cdot i^{k_s-k_t}\prod_{j=1}^r\omega_j^{v_j+u_j-w_j}\overline{\omega}_j^{v_j-u_j+w_j}\in\mathbb{Z}[i]
\end{eqnarray*}
we have that the distance between $\xi_s$ and $\xi_t$ is integral for every $1\le s,t\le 2\tau(R)$. Additionally we can add the center of the circle to this point set to obtain an integral point set of cardinality $2\cdot\tau(R)+1$ having rational coordinates.
After a suitable rotation we can achieve integral coordinates.

So let us have an example. We choose $R=5\cdot 13=65$ and successively obtain
\begin{eqnarray*}
  \omega_1=2+i,\quad
  \omega_2=3+2i,\\
  \eta_{1}  = 65i,    \quad
  \eta_{2}  = 65,     \quad
  \eta_{3}  = -52+39i,\quad
  \eta_{4}  = 39+52i, \\
  \eta_{5} = -60+25i, \quad
  \eta_{6} = 25+60i,  \quad
  \eta_{7} = -56-33i, \quad
  \eta_{8} = -33+56i, \\
  \xi_{1}=-65,        \quad
  \xi_{2}=65,         \quad
  \xi_{3}=\frac{91}{5}-\frac{312}{5}i,\quad
  \xi_{4}=-\frac{91}{5}+\frac{312}{5}i,\\
  \xi_{5}=\frac{595}{13}-\frac{600}{13}i,\quad
  \xi_{6}=-\frac{595}{13}+\frac{600}{13}i,\quad
  \xi_{7}=\frac{2047}{65}+\frac{3696}{65}i,\quad
  \xi_{8}=-\frac{2047}{65}-\frac{3696}{65}i.
\end{eqnarray*}
After adding the origin $(0,0)^T$ and applying a suitable rotation and translation we obtain the maximal integral point set
\begin{eqnarray*}
  \mathcal{P}&=&\left[
  \begin{pmatrix} 0\\0\end{pmatrix},
  \begin{pmatrix} 0\\-32\end{pmatrix},
  \begin{pmatrix} -30\\40\end{pmatrix},
  \begin{pmatrix} -30\\-72\end{pmatrix},
  \begin{pmatrix} -63\\-16\end{pmatrix},
  \begin{pmatrix} -96\\40\end{pmatrix},
  \begin{pmatrix} -96\\-72\end{pmatrix},
  \begin{pmatrix} -126\\0\end{pmatrix},
  \begin{pmatrix} -126\\-32\end{pmatrix}
  \right]
\end{eqnarray*}
in minimum coordinate representation.

\begin{Construction}
  \label{circle_construction}
  For a given $R$ which has only prime factors $p$ fulfilling $p\equiv 1\pmod 4$ there exists an integral point set $\mbox{circle}(R)$
  consisting of $2\cdot\tau(R)$ points on a circle of radius $R$ together with its center, where $\tau(R)$ denotes the number
  of divisors of $R$.
\end{Construction}

From the above it is easy to deduce that the $2\tau(R)$ points on the circle all have pairwise even distances and that the diameter of this point set is given by $2R$. Using this we can give another construction.

\begin{Construction}
  \label{circle_construction_tilde}
  For a given $R$ which has only prime factors $p$ fulfilling $p\equiv 1\pmod 4$ there exists an integral point set
  $\widetilde{\mbox{circle}}(R)$ consisting of $2\cdot\tau(R)$ points on a circle of radius $\frac{R}{2}$, where $\tau(R)$ denotes
  the number of divisors of $R$.
\end{Construction}

\begin{Conjecture}
  The plane integral point sets given by Construction \ref{circle_construction} and Construction \ref{circle_construction_tilde} are
  maximal.
\end{Conjecture}

We can generalize the idea of Construction \ref{circle_construction_tilde} in some way. Let $t$ be an arbitrary integer, $R$ be a
integer having only prime factor fulfilling $p\equiv 1\pmod 4$, and $\mathcal{P}(R)$ be the integral point set given by Construction
\ref{circle_construction} with radius $R$. By $\mathcal{P}(R,t)$ we denote the point set which arises from $\mathcal{P}(R)$ by scaling the point set with a factor $\frac{1}{t}$, this means dividing all distances by $t$. Thus $\mathcal{P}(R,t)$ is a point set
with pairwise rational distances and rational coordinates. With this we can construct a graph $\mathcal{G}$ containing the points
of $\mathcal{P}(R,t)$ as its vertices. Two vertices of $\mathcal{G}$ are connected by an edge, if and only if the corresponding points
have an integral distance in $\mathcal{P}(R,t)$. The maximal cliques $\mathcal{C}$ of $\mathcal{G}$ correspond to integral point sets
$\mathcal{P}(R,t,\mathcal{C})$.

\begin{Construction}
\label{circle_construction_scaling}
  For a given $R$ which has only prime factors $p$ fulfilling $p\equiv 1\pmod 4$ and a given integer $t$ there exist
  integral point sets $\mbox{circle}(R,t,\mathcal{C})$ consisting of points on a circle of radius $\frac{R}{t}$,
  where $\mathcal{C}$ is a maximal clique of the above described graph. As an abbreviation we use $\mbox{circle}(R,t)$
  instead of $\mbox{circle}(R,t,\mathcal{C})$.
\end{Construction}

\begin{Conjecture}
  For $t=8$ Construction \ref{circle_construction_scaling} gives maximal integral point sets of cardinality $\tau(R)$.
\end{Conjecture}

\section{Maximal integral point sets over $\mathbb{Z}^2$ with further conditions}
\label{sec_further_conditions}

\begin{table}[htp]
  \begin{center}
    \begin{tabular}{|r|r|l||r|r|l|}
      \hline
      $k$ & $d_M(k,2)$ & construction & $k$ & $d_M(k,2)$ & construction \\
      \hline
         3 &                                     $=2066$ & $\Delta(2066,1803,505)$ &
        26 &                                $\le 112895$ & $\mbox{decompose}\left(2^6\cdot 3\cdot 7,5\right)$\\
         4 &                                        $=5$ & $\mathcal{P}_1(3,4)=\widetilde{\mbox{circle}}(5)$ &
        27 &                                  $\le 2590$ & $\mbox{decompose}\left(2^3\cdot 3^2\right)$\\
         5 &                                        $=8$ & $\mathcal{P}_2(3,4)=\mbox{crab}(3,4)$ &
        28 &          $\overset{\star}{\lesssim} 203125$ & $\widetilde{\mbox{circle}}\left(5^{6}\cdot 13\right)$\\
           &                                             & $=\mbox{crab}(4,3)$ &
        29 &                                  $\le 1798$ & $\mbox{decompose}\left(2^2\cdot 3\cdot 5\right)$\\
         6 &                                       $=25$ & $\widetilde{\mbox{circle}}\left(5^2\right)$ &
        30 &                                $\le 105625$ & $\widetilde{\mbox{circle}}\left(5^{4}\cdot 13^2\right)$\\
         7 &                                       $=30$ & $\mbox{crab}\left(8,6,15\right)$ &
        31 &          $\overset{\star}{\lesssim} 211250$ & $\mbox{circle}\left(5^{4}\cdot 13^2\right)$\\
         8 &                                       $=65$ & $\widetilde{\mbox{circle}}\left(5\cdot 13\right)$ &
        32 &                                 $\le 27625$ & $\widetilde{\mbox{circle}}\left(5^{3}\cdot 13\cdot 17\right)$\\
         9 &                                      $=130$ & $\mbox{circle}\left(5\cdot 13\right)$ &
        33 &                                 $\le 55250$ & $\mbox{circle}\left(5^{3}\cdot 13\cdot 17\right)$\\
        10 &                                   $\le 625$ & $\widetilde{\mbox{circle}}\left(5^4\right)$ &
        34 &                                $\le 142295$ & $\mbox{decompose}\left(2^3\cdot 3\cdot 7\cdot 11,5\right)$\\
        11 &                                       $=70$ & $\mbox{decompose}\left(2^2\cdot 3\right)$ &
        35 &                                 $\le 18430$ & $\mbox{decompose}\left(2^6\cdot 3\right)$\\
        12 &                                     $= 325$ & $\widetilde{\mbox{circle}}\left(5^2\cdot 13\right)$ &
        36 &                                 $\le 40625$ & $\widetilde{\mbox{circle}}\left(5^{5}\cdot 13\right)$\\
        13 &                                   $\le 650$ & $\mbox{circle}\left(5^2\cdot 13\right)$ &
        37 &                                 $\le 10366$ & $\mbox{decompose}\left(2^4\cdot 3^2\right)$\\
        14 &                                 $\le 15625$ & $\widetilde{\mbox{circle}}\left(5^6\right)$ &
        38 &          $\overset{\star}{\lesssim} 571535$ & $\mbox{decompose}\left(2^4\cdot 3^3\cdot 7,5\right)$\\
        15 &                                  $\le 8190$ & $\mbox{decompose}\left(2^7\right)$ &
        39 &         $\overset{\star}{\lesssim} 4816895$ & $\mbox{decompose}\left(2^{9}\cdot 3\cdot 7,5\right)$\\
        16 &                                  $\le 1105$ & $\widetilde{\mbox{circle}}\left(5\cdot 13\cdot 17\right)$ &
        40 &                                $\le 138125$ & $\widetilde{\mbox{circle}}\left(5^{4}\cdot 13\cdot 17\right)$\\
        17 &                                     $= 286$ & $\mbox{decompose}\left(2^3\cdot 3\right)$ &
        41 &                                 $\le 73726$ & $\mbox{decompose}\left(2^7\cdot 3\right)$\\
        18 &                                  $\le 4225$ & $\widetilde{\mbox{circle}}\left(5^2\cdot 13^2\right)$ &
        42 &          $\overset{\star}{\lesssim} 677375$ & $\mbox{decompose}\left(2^{6}\cdot 3^2\cdot 7,5\right)$\\
        19 &                                  $\le 8450$ & $\mbox{circle}\left(5^2\cdot 13^2\right)$ &
        43 &         $\overset{\star}{\lesssim} 4573799$ & $\mbox{decompose}\left(2^3\cdot 3^2\cdot 5\cdot 7\cdot 11,17\right)$\\
        20 &                                  $\le 8125$ & $\widetilde{\mbox{circle}}\left(5^4\cdot 13\right)$ &
        44 &         $\overset{\star}{\lesssim} 6614998$ & $\mbox{decompose}\left(2^4\cdot 3^2\cdot 5^2\cdot 7,13\right)$\\
        21 &                                 $\le 16250$ & $\mbox{circle}\left(5^4\cdot 13\right)$ &
        45 &         $\overset{\star}{\lesssim} 7001315$ & $\mbox{decompose}\left(2^3\cdot 3^2\cdot 7^2,5\right)$\\
        22 &                                 $\le 53360$ & $\mbox{decompose}\left(2^2\cdot 3\cdot 7\cdot 11,5\right)$ &
        46 &        $\overset{\star}{\lesssim} 64833614$ & $\mbox{decompose}\left(2^2\cdot 3^4\cdot 5\cdot 7\cdot 11,17\right)$\\
        23 &                                  $\le 1150$ & $\mbox{decompose}\left(2^4\cdot 3\right)$ &
        47 &                                  $\le 7198$ & $\mbox{decompose}\left(2^{3}\cdot 3\cdot 5\right)$\\
        24 &                                  $\le 5525$ & $\widetilde{\mbox{circle}}\left(5^{2}\cdot 13\cdot 17\right)$ &
        48 &          $\overset{\star}{\lesssim} 160225$ & $\widetilde{\mbox{circle}}\left(5^{2}\cdot 13\cdot 17\cdot 29\right)$\\
        25 &                                 $\le 11050$ & $\mbox{circle}\left(5^{2}\cdot 13\cdot 17\right)$ &
        49 &          $\overset{\star}{\lesssim} 320450$ & $\mbox{circle}\left(5^{2}\cdot 13\cdot 17\cdot 29\right)$\\
           &                                             &                                                           &
        50 &         $\overset{\star}{\lesssim} 4064255$ & $\mbox{decompose}\left(2^{7}\cdot 3^2\cdot 7,5 \right)$\\
      \hline
    \end{tabular}
    \caption{Best known constructions for maximal integral point sets over $\mathbb{Z}^2$ in arbitrary position.}
    \label{table_min_arbitrary_position}
  \end{center}
\end{table}

In Table \ref{table_min_arbitrary_position} we have summarized the constructions yielding the smallest diameter of a maximal integral point set over $\mathbb{Z}^2$. Some of the values $d_M(k,2)$ could be determined exactly by an exhaustive search, but for most values of $k$ we only have upper bounds (and $301$ as lower bound). In some cases, denoted by $\overset{\star}{\lesssim}$, we were not able to check the maximality of the constructed point sets, since their diameter was too large. Looking at Table \ref{table_min_arbitrary_position} we observe, that the constructions of crabs (Construction \ref{line_construction} and Construction \ref{semi_crab_construction}) are very dominating. The resulting point sets contain $n-2$ and $n-1$ collinear points out of $n$ points, respectively. So it may be interesting to study maximal integral point sets over $\mathbb{Z}^2$, where no three points are collinear. We also say, that a point set is in semi-general position, if no three points are collinear. By $\overline{d_M}(k,2)$ we denote the minimum possible diameter of these point sets. We can check for this further condition, that no three points are collinear, by applying Lemma \ref{lemma_collinear}. Using the methods and algorithms described in this article, we were able to obtain some exact values and some upper bounds for $\overline{d_M}(k,2)$. The results are summarized in Table \ref{table_min_semi_general_position}. We would like to remark that we additionally have the lower bounds $\overline{d_M}(k,2)\ge 5525$ for $k\in\{11,13,14,15,17\}$ and $\overline{d_M}(k,2)\ge 10001$ for $k\ge 19$, $k\neq 20,24$.

\begin{table}[htp]
  \begin{center}
    \begin{tabular}{|r|r|l||r|r|l|}
      \hline
      $k$ & $\overline{d_M}(k,2)$ & construction & $k$ & $\overline{d_M}(k,2)$ & construction \\
      \hline
         3 &                                       $=2066$ & $\scriptstyle{\Delta(2066,1803,505)}$ &
        27 &             $\overset{\star}{\lesssim}305218$ & $\scriptstyle{\mbox{circle}\left(5^{2}\cdot 13^2\cdot 17^2,8\right)}$ \\
         4 &                                          $=5$ & $\scriptstyle{\mathcal{P}_1(3,4)=\widetilde{\mbox{circle}}(5)}$ &
        28 &            $\overset{\star}{\lesssim} 203125$ & $\scriptstyle{\widetilde{\mbox{circle}}\left(5^{6}\cdot 13\right)}$\\
         5 &                                        $=120$ & see Figure \ref{fig:semi_general_5_smallest}&
        29 &$\overset{\star}{\lesssim}9311389618298531250$ & $\scriptstyle{\mbox{circle}\left(5^{28},8\right)}$ \\
         6 &                                         $=25$ & $\scriptstyle{\widetilde{\mbox{circle}}\left(5^2\right)}$ &
        30 &                                  $\le 105625$ & $\scriptstyle{\widetilde{\mbox{circle}}\left(5^{4}\cdot 13^2\right)}$\\
         7 &                                        $=925$ &  see Figure \ref{fig:semi_general_7_smallest}&
        31 & $\overset{\star}{\lesssim}232784740457463281250$ & $\scriptstyle{\mbox{circle}\left(5^{30},8\right)}$ \\
         8 &                                         $=65$ & $\scriptstyle{\widetilde{\mbox{circle}}\left(5\cdot 13\right)}$ &
        32 &                                   $\le 27625$ & 
             $\scriptstyle{\widetilde{\mbox{circle}}\left(5^{3}\cdot 13\cdot 17\right)}$\\
         9 &                                       $=1045$ & see Figure \ref{fig:semi_general_9_smallest}&
        33 &          $\overset{\star}{\lesssim}412343750$ & $\scriptstyle{\mbox{circle}\left(5^{10}\cdot 13^2,8\right)}$ \\
        10 &                                        $=625$ & $\scriptstyle{\widetilde{\mbox{circle}}\left(5^4\right)}$ &
        34 &      $\overset{\star}{\lesssim} 152587890625$ & $\scriptstyle{\widetilde{\mbox{circle}}\left(5^{16}\right)}$\\
        11 &            $\overset{\star}{\lesssim}2434375$ & $\scriptstyle{\mbox{circle}\left(5^{10},8\right)}$ &
        35 &          $\overset{\star}{\lesssim}111562500$ & $\scriptstyle{\mbox{circle}\left(5^{6}\cdot 13^4,8\right)}$ \\
        12 &                                        $=325$ & $\scriptstyle{\widetilde{\mbox{circle}}\left(5^2\cdot 13\right)}$ &
        36 &                                   $\le 71825$ & 
             $\scriptstyle{\widetilde{\mbox{circle}}\left(5^{2}\cdot 13^2\cdot 17\right)}$\\
        13 &           $\overset{\star}{\lesssim}60859375$ & $\scriptstyle{\mbox{circle}\left(5^{12},8\right)}$ &
        37 & $\scriptstyle{\overset{\star}{\lesssim}3637261569647863769531250}$ &
             $\scriptstyle{\mbox{circle}\left(5^{36},8\right)}$ \\
        14 &                                   $\le 15625$ & $\scriptstyle{\widetilde{\mbox{circle}}\left(5^6\right)}$ &
        38 &     $\overset{\star}{\lesssim} 3814697265625$ & $\scriptstyle{\widetilde{\mbox{circle}}\left(5^{18}\right)}$\\
        15 &                                   $\le 26390$ & $\scriptstyle{\mbox{circle}\left(5^{4}\cdot 13^2,8\right)}$ &
        39 &        $\overset{\star}{\lesssim}10314771205$ & $\scriptstyle{\mbox{circle}\left(5^{12}\cdot 13^2,8\right)}$ \\
        16 &                                       $=1105$ & $\scriptstyle{\widetilde{\mbox{circle}}\left(5\cdot 13\cdot 17\right)}$ &
        40 &                                  $\le 138125$ & 
             $\scriptstyle{\widetilde{\mbox{circle}}\left(5^{4}\cdot 13\cdot 17\right)}$\\
        17 &        $\scriptstyle{\overset{\star}{\lesssim}38037109375}$ & $\scriptstyle{\mbox{circle}\left(5^{16},8\right)}$ &
        41 & $\scriptstyle{\overset{\star}{\lesssim}}$\tiny{2273288481029914855957031250} &
             $\scriptstyle{\mbox{circle}\left(5^{40},8\right)}$ \\
        18 &                                       $=4225$ & $\scriptstyle{\widetilde{\mbox{circle}}\left(5^2\cdot 13^2\right)}$ &
        42 &           $\overset{\star}{\lesssim} 2640625$ & $\scriptstyle{\widetilde{\mbox{circle}}\left(5^{6}\cdot 13^2\right)}$\\
        19 &       $\scriptstyle{\overset{\star}{\lesssim}950927734375}$ & $\scriptstyle{\mbox{circle}\left(5^{18},8\right)}$ &
        43 & $\scriptstyle{\overset{\star}{\lesssim}}$\tiny{56832212025747871398925781250} &
             $\scriptstyle{\mbox{circle}\left(5^{42},8\right)}$ \\
        20 &                                       $=8125$ & $\scriptstyle{\widetilde{\mbox{circle}}\left(5^4\cdot 13\right)}$ &
        44 &         $\overset{\star}{\lesssim} 126953125$ & $\scriptstyle{\widetilde{\mbox{circle}}\left(5^{10}\cdot 13\right)}$\\
        21 &             $\overset{\star}{\lesssim}659750$ & $\scriptstyle{\mbox{circle}\left(5^{6}\cdot 13^2,8\right)}$ &
        45 &            $\overset{\star}{\lesssim}7630450$ & $\scriptstyle{\mbox{circle}\left(5^{4}\cdot 13^2\cdot 17^2,8\right)}$ \\
        22 &           $\overset{\star}{\lesssim} 9765625$ & $\scriptstyle{\widetilde{\mbox{circle}}\left(5^{10}\right)}$ &
        46 &  $\overset{\star}{\lesssim} 2384185791015625$ & $\scriptstyle{\widetilde{\mbox{circle}}\left(5^{22}\right)}$\\
        23 & \!\!$\scriptstyle{\overset{\star}{\lesssim}}$\tiny{595928935571106} &
             $\scriptstyle{\mbox{circle}\left(5^{22},8\right)}$ &
        47 & $\scriptstyle{\overset{\star}{\lesssim}}$\tiny{35520132516092419624328613281250} &
             $\scriptstyle{\mbox{circle}\left(5^{46},8\right)}$ \\
        24 &                                       $=5525$ & 
             $\scriptstyle{\widetilde{\mbox{circle}}\left(5^{2}\cdot 13\cdot 17\right)}$ &
        48 &            $\overset{\star}{\lesssim} 160225$ & 
             $\scriptstyle{\widetilde{\mbox{circle}}\left(5^{2}\cdot 13\cdot 17\cdot 29\right)}$\\
        25 &            $\overset{\star}{\lesssim}4462500$ & $\scriptstyle{\mbox{circle}\left(5^{4}\cdot 13^4,8\right)}$ &
        49 &        $\overset{\star}{\lesssim}18854062500$ & $\scriptstyle{\mbox{circle}\left(5^{6}\cdot 13^6,8\right)}$ \\
        26 &         $\overset{\star}{\lesssim} 244140625$ & $\scriptstyle{\widetilde{\mbox{circle}}\left(5^{12}\right)}$ &
        50 &          $\overset{\star}{\lesssim} 17850625$ & $\scriptstyle{\widetilde{\mbox{circle}}\left(5^{4}\cdot 13^4 \right)}$\\
      \hline
    \end{tabular}
    \caption{Best known constructions for maximal integral point sets over $\mathbb{Z}^2$ in semi-general position.}
    \label{table_min_semi_general_position}
  \end{center}
\end{table}

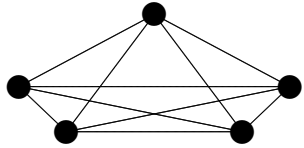
\begin{figure}[htp]
  \begin{center}
    \setlength{\unitlength}{0.3mm}
    \begin{picture}(120,52)
      \put(0,20){\circle*{10}}
      \put(120,20){\circle*{10}}
      \put(99,0){\circle*{10}}
      \put(60,52){\circle*{10}}
      \put(21,0){\circle*{10}}
      \qbezier(0,20)(60,20)(120,20)
      \qbezier(0,20)(49.5,10)(99,0)
      \qbezier(0,20)(30,36)(60,52)
      \qbezier(0,20)(10.5,10,)(21,0)
      \qbezier(120,20)(109.5,10)(99,0)
      \qbezier(120,20)(90,36)(60,52)
      \qbezier(120,20)(70.5,10)(21,0)
      \qbezier(99,0)(79.5,26)(60,52)
      \qbezier(99,0)(60,0)(21,0)
      \qbezier(60,52)(40.5,26)(21,0)
    \end{picture}
    \caption{The smallest maximal integral point set of cardinality $5$ in semi-general position.}
    \label{fig:semi_general_5_smallest}
  \end{center}
\end{figure}

We would like to have a closer look on the smallest known examples of maximal integral point sets in semi-general position consisting of an odd number of points. For cardinality $5$ the two smallest point sets with respect to the diameter are given in minimum coordinate representation by
$$
  \left[
    \begin{pmatrix}0\\0\end{pmatrix},\\
    \begin{pmatrix}0\\-78\end{pmatrix},\\
    \begin{pmatrix}-20\\21\end{pmatrix},\\
    \begin{pmatrix}-20\\-99\end{pmatrix},\\
    \begin{pmatrix}-52,-39\end{pmatrix}
  \right]
  \text{ and }
  \left[
    \begin{pmatrix}0\\0\end{pmatrix},\\
    \begin{pmatrix}0\\-80\end{pmatrix},\\
    \begin{pmatrix}-45\\28\end{pmatrix},\\
    \begin{pmatrix}-45\\-108\end{pmatrix},\\
    \begin{pmatrix}-96\\-40\end{pmatrix}
  \right],
$$
see Figure \ref{fig:semi_general_5_smallest} for a drawing of the first point set. Both point sets consist of four point on a circle $\mathcal{C}$ of radii $\frac{29\cdot 101}{40}$ and $\frac{13\cdot 53}{10}$, respectively. In each case the fifth point does not lie on this circle $\mathcal{C}$, but the line through this point and the center of $\mathcal{C}$ is a symmetry axis of the point set.

\begin{figure}[htp]
  \begin{center}
    \setlength{\unitlength}{0.075mm}
    \begin{picture}(624,880)
      \put(99,0){\circle*{30}}
      \put(384,880){\circle*{30}}
      \put(99,880){\circle*{30}}
      \put(624,700){\circle*{30}}
      \put(624,180){\circle*{30}}
      \put(0,440){\circle*{30}}
      \put(384,0){\circle*{30}}
      \qbezier(99,0)(241.5,440)(384,880)
      \qbezier(99,0)(99,440)(99,880)
      \qbezier(99,0)(361.5,350)(624,700)
      \qbezier(99,0)(361.5,90)(624,180)
      \qbezier(99,0)(49.5,220)(0,440)
      \qbezier(99,0)(241.5,0)(384,0)
      \qbezier(384,880)(241.5,880)(99,880)
      \qbezier(384,880)(504,790)(624,700)
      \qbezier(384,880)(504,530)(624,180)
      \qbezier(384,880)(192,660)(0,440)
      \qbezier(384,880)(384,440)(384,0)
      \qbezier(99,880)(361.5,790)(624,700)
      \qbezier(99,880)(361.5,530)(624,180)
      \qbezier(99,880)(49.5,660)(0,440)
      \qbezier(99,880)(241.5,440)(384,0)
      \qbezier(624,700)(624,440)(624,180)
      \qbezier(624,700)(312,570)(0,440)
      \qbezier(624,700)(504,350)(384,0)
      \qbezier(624,180)(312,310)(0,440)
      \qbezier(624,180)(504,90)(384,0)
      \qbezier(0,440)(192,220)(384,0)
    \end{picture}
    \caption{The smallest maximal integral point set of cardinality $7$ in semi-general position.}
    \label{fig:semi_general_7_smallest}
  \end{center}
\end{figure}
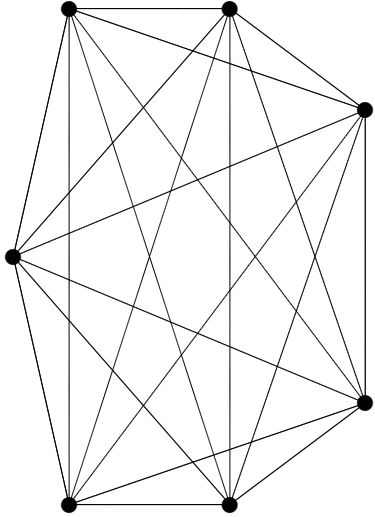

\noindent
For cardinality $7$ the two smallest examples are given by
$$
  \left[
    \begin{pmatrix}0\\0\end{pmatrix},\\
    \begin{pmatrix}0\\-285\end{pmatrix},\\
    \begin{pmatrix}-180\\240\end{pmatrix},\\
    \begin{pmatrix}-440\\-384\end{pmatrix},\\
    \begin{pmatrix}-700\\240\end{pmatrix},\\
    \begin{pmatrix}-880\\0\end{pmatrix},\\
    \begin{pmatrix}-880\\-285\end{pmatrix}
  \right]
$$
and
$$
  \left[
    \begin{pmatrix}0\\0\end{pmatrix},\\
    \begin{pmatrix}0\\-855\end{pmatrix},\\
    \begin{pmatrix}-540\\720\end{pmatrix},\\
    \begin{pmatrix}-1320\\-1152\end{pmatrix},\\
    \begin{pmatrix}-2100\\720\end{pmatrix},\\
    \begin{pmatrix}-2640\\0\end{pmatrix},\\
    \begin{pmatrix}-2640\\-855\end{pmatrix}
  \right],
$$
see Figure \ref{fig:semi_general_7_smallest} for a graphical representation of the first example. The geometric shape of the corresponding two point sets is similar to the case of cardinality $5$. In each case $6$ points are situated on a circle $\mathcal{C}$ of radii $\frac{5^2\cdot 37}{2}$ and $\frac{3\cdot 5^2\cdot 37}{2}$, respectively. Again we have the symmetry axis through the seventh point and the center of $\mathcal{C}$.

\begin{figure}[htp]
  \begin{center}
    \setlength{\unitlength}{0.075mm}
    \begin{picture}(1014,900)
      \put(255,64){\circle*{30}}
      \put(882,900){\circle*{30}}
      \put(969,784){\circle*{30}}
      \put(1014,676){\circle*{30}}
      \put(132,900){\circle*{30}}
      \put(45,784){\circle*{30}}
      \put(0,676){\circle*{30}}
      \put(759,64){\circle*{30}}
      \put(507,0){\circle*{30}}
      \qbezier(255,64)(568.5,482)(882,900)
      \qbezier(255,64)(612,424)(969,784)
      \qbezier(255,64)(634.5,370)(1014,676)
      \qbezier(255,64)(193.5,482)(132,900)
      \qbezier(255,64)(150,424)(45,784)
      \qbezier(255,64)(127.5,370)(0,676)
      \qbezier(255,64)(507,64)(759,64)
      \qbezier(255,64)(381,32)(507,0)
      \qbezier(882,900)(925.5,842)(969,784)
      \qbezier(882,900)(948,788)(1014,676)
      \qbezier(882,900)(507,900)(132,900)
      \qbezier(882,900)(463.5,842)(45,784)
      \qbezier(882,900)(441,788)(0,676)
      \qbezier(882,900)(820.5,482)(759,64)
      \qbezier(882,900)(694.5,450)(507,0)
      \qbezier(969,784)(991.5,730)(1014,676)
      \qbezier(969,784)(550.5,842)(132,900)
      \qbezier(969,784)(507,784)(45,784)
      \qbezier(969,784)(484.5,730)(0,676)
      \qbezier(969,784)(864,424)(759,64)
      \qbezier(969,784)(738,392)(507,0)
      \qbezier(1014,676)(573,788)(132,900)
      \qbezier(1014,676)(529.5,730)(45,784)
      \qbezier(1014,676)(507,676)(0,676)
      \qbezier(1014,676)(886.5,370)(759,64)
      \qbezier(1014,676)(760.5,338)(507,0)
      \qbezier(132,900)(88.5,842)(45,784)
      \qbezier(132,900)(66,788)(0,676)
      \qbezier(132,900)(445.5,482)(759,64)
      \qbezier(132,900)(319.5,450)(507,0)
      \qbezier(45,784)(22.4,730)(0,676)
      \qbezier(45,784)(402,424)(759,64)
      \qbezier(45,784)(276,392)(507,0)
      \qbezier(0,676)(379.5,370)(759,64)
      \qbezier(0,676)(253.5,338)(507,0)
      \qbezier(759,64)(633,32)(507,0)
    \end{picture}
    \caption{The smallest maximal integral point set of cardinality $9$ in semi-general position.}
    \label{fig:semi_general_9_smallest}
  \end{center}
\end{figure}
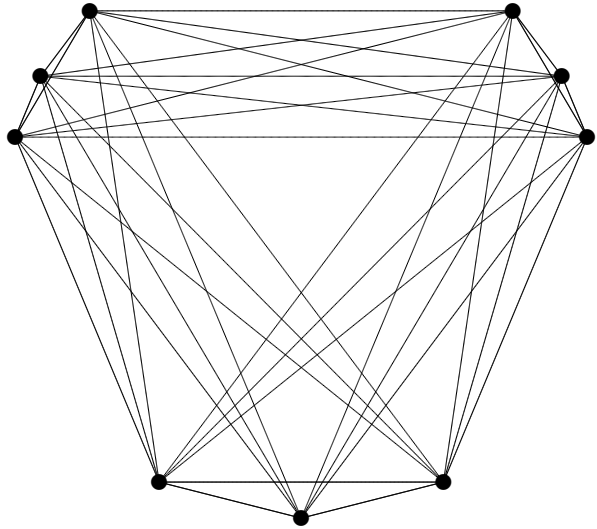

\noindent
For cardinality $9$ the two smallest examples are given by
$$
  \left[
    \begin{pmatrix}0\\0\end{pmatrix},\\
    \begin{pmatrix}0\\-504\end{pmatrix},\\
    \begin{pmatrix}-64\\-252\end{pmatrix},\\
    \begin{pmatrix}612\\255\end{pmatrix},\\
    \begin{pmatrix}612\\-759\end{pmatrix},\\
    \begin{pmatrix}720\\210\end{pmatrix},\\
    \begin{pmatrix}720\\-714\end{pmatrix},\\
    \begin{pmatrix}836\\123\end{pmatrix},\\
    \begin{pmatrix}836\\-627\end{pmatrix}
  \right]
$$
and
$$
  \left[
    \begin{pmatrix}0\\0\end{pmatrix},\\
    \begin{pmatrix}0\\-672\end{pmatrix},\\
    \begin{pmatrix}-123\\164\end{pmatrix},\\
    \begin{pmatrix}-123\\-836\end{pmatrix},\\
    \begin{pmatrix}-816\\340\end{pmatrix},\\
    \begin{pmatrix}-816\\-1012\end{pmatrix},\\
    \begin{pmatrix}-960\\280\end{pmatrix},\\
    \begin{pmatrix}-960\\-952\end{pmatrix},\\
    \begin{pmatrix}-1323\\-336\end{pmatrix}
  \right],
$$
see Figure \ref{fig:semi_general_9_smallest} for a graphical representation of the first example. Here in both examples all nine points are situated of circles of radii $\frac{5^2\cdot 13^2}{8}$ and $\frac{5^2\cdot 13^2}{6}$, respectively. They both can be
obtained using Construction \ref{circle_construction_scaling}.

Now we observe that the constructions based on circles, Construction \ref{circle_construction}, Construction \ref{circle_construction_tilde}, and Construction \ref{circle_construction_scaling}, are very dominating in this context. The next natural step is to also forbid four points on a circle. If no three points are on a line and no four points on a circle we speak of general position. By $\dot{d}_M(k,2)$ we denote the minimum possible diameter of a maximal plane integral point set in general position over $\mathbb{Z}^2$. Without the maximality condition these point sets are also known as $k_2$-cluster \cite{0676.52006}. As we cannot apply our most successful constructions based on crabs and circles in this case, examples are scarce. For the check whether four points are situated on a circle we have a well known criterion similar to Lemma \ref{lemma_collinear}:
\begin{Lemma}
  \label{lemma_four_on_circle}
  Four points $(x_1,y_1)$, $(x_2,y_2)$, $(x_3,y_3)$, $(x_4,y_4)$ in $\mathbb{R}^2$ are situated on a circle if and only if
  $$
    \left|
    \begin{array}{cccc}
      x_1 & y_1 & x_1^2+y_1^2 & 1 \\
      x_2 & y_2 & x_2^2+y_2^2 & 1 \\
      x_3 & y_3 & x_3^2+y_3^2 & 1 \\
      x_4 & y_4 & x_4^2+y_4^2 & 1
    \end{array}
    \right|=0
  $$
  holds.
\end{Lemma}

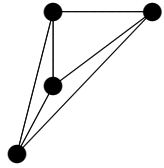
\begin{figure}[htp]
  \begin{center}
    \setlength{\unitlength}{0.3mm}
    \begin{picture}(60,63)
      \put(0,0){\circle*{8}}
      \put(60,63){\circle*{8}}
      \put(16,63){\circle*{8}}
      \put(16,30){\circle*{8}}
      \qbezier(0,0)(30,31.5)(60,63)
      \qbezier(0,0)(8,31.5)(16,63)
      \qbezier(0,0)(8,15)(16,30)
      \qbezier(60,63)(38,63)(16,63)
      \qbezier(60,63)(38,46.5)(16,30)
      \qbezier(16,63)(16,46.5)(16,30)
    \end{picture}
    \caption{The smallest maximal integral point set of cardinality $4$ in general position.}
    \label{fig:general_4_smallest}
  \end{center}
\end{figure}

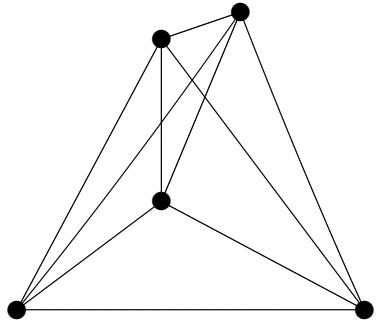
\begin{figure}[htp]
  \begin{center}
    \setlength{\unitlength}{0.3mm}
    \begin{picture}(154,132)
      \put(0,0){\circle*{8}}
      \put(99,132){\circle*{8}}
      \put(154,0){\circle*{8}}
      \put(64,120){\circle*{8}}
      \put(64,48){\circle*{8}}
      \qbezier(0,0)(49.5,61)(99,132)
      \qbezier(0,0)(77,0)(154,0)
      \qbezier(0,0)(32,60)(64,120)
      \qbezier(0,0)(32,24)(64,48)
      \qbezier(99,132)(126.5,61)(154,0)
      \qbezier(99,132)(81.5,126)(64,120)
      \qbezier(99,132)(81.5,90)(64,48)
      \qbezier(154,0)(109,60)(64,120)
      \qbezier(154,0)(109,24)(64,48)
      \qbezier(64,120)(64,84)(64,48)
    \end{picture}
    \caption{The smallest maximal integral point set of cardinality $5$ in general position.}
    \label{fig:general_5_smallest}
  \end{center}
\end{figure}

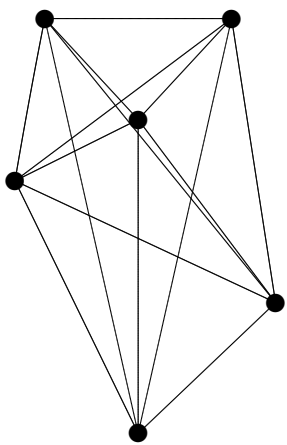
\begin{figure}[htp]
  \begin{center}
    \setlength{\unitlength}{0.03mm}
    \begin{picture}(1155,1840)
      \put(546,0){\circle*{80}}
      \put(960,1840){\circle*{80}}
      \put(132,1840){\circle*{80}}
      \put(546,1392){\circle*{80}}
      \put(0,1120){\circle*{80}}
      \put(1155,580){\circle*{80}}
      \qbezier(546,0)(753,920)(960,1840)
      \qbezier(546,0)(339,920)(132,1840)
      \qbezier(546,0)(546,696)(546,1392)
      \qbezier(546,0)(273,560)(0,1120)
      \qbezier(546,0)(850.5,290)(1155,580)
      \qbezier(960,1840)(546,1840)(132,1840)
      \qbezier(960,1840)(753,1616)(546,1392)
      \qbezier(960,1840)(480,1480)(0,1120)
      \qbezier(960,1840)(1057.5,1210)(1155,580)
      \qbezier(132,1840)(339,1616)(546,1392)
      \qbezier(132,1840)(61,1480)(0,1120)
      \qbezier(132,1840)(643.5,1210)(1155,580)
      \qbezier(546,1392)(273,1256)(0,1120)
      \qbezier(546,1392)(850.5,986)(1155,580)
      \qbezier(0,1120)(577.5,850)(1155,580)
    \end{picture}
    \caption{The smallest maximal integral point set of cardinality $6$ in general position.}
    \label{fig:general_6_smallest}
  \end{center}
\end{figure}

In Table \ref{table_min_general_position} we have summarized our knowledge on $\dot{d}_M(k,2)$. For the lower bound $\dot{d}_M(7,2)>599000$ we refer to \cite{herontriangles}. Whether $\dot{d}_M(7,2)$ is finite (even if we drop the maximality condition) is an open problem, see \cite{UPIN,0676.52006}. If we drop the maximality condition and the condition on the integrality of the coordinates (in other words characteristic one), then very recently two such examples were found, see \cite{paper_kreisel}. The smallest example for $k=6$ is indeed the smallest integral point set of characteristic one in general position with cardinality $6$. We would also like to give the coordinates for the second smallest examples. For cardinality $4$ we have
$$
  \left[
  \begin{pmatrix}0\\0\end{pmatrix},
  \begin{pmatrix}0\\-69\end{pmatrix},
  \begin{pmatrix}-20\\-21\end{pmatrix}, \\
  \begin{pmatrix}-92\\0\end{pmatrix}
  \right],
$$
for cardinality $5$ we have
$$
  \left[
  \begin{pmatrix}0\\0\end{pmatrix},
  \begin{pmatrix}0\\-153\end{pmatrix},
  \begin{pmatrix}-60\\144\end{pmatrix},
  \begin{pmatrix}-140\\-48\end{pmatrix},
  \begin{pmatrix}-176\\57\end{pmatrix}
  \right],
$$
and for cardinality $6$ we have
$$
  \left[
  \begin{pmatrix}0\\0\end{pmatrix},
  \begin{pmatrix}-135\\-1008\end{pmatrix},
  \begin{pmatrix}420\\1008\end{pmatrix},
  \begin{pmatrix}735\\-392\end{pmatrix},
  \begin{pmatrix}1155\\616\end{pmatrix},
  \begin{pmatrix}1290\\1624\end{pmatrix}
  \right].
$$

\begin{table}[htp]
  \begin{center}
    \begin{tabular}{|r|r|l|}
      \hline
      $k=|\mathcal{P}|$ & $\dot{d}_M(k,2)$ & construction \\
      \hline
         3 &                                     $=2066$ & $\Delta(2066,1803,505)$\\
         4 &                                       $=87$ & $\left[\begin{pmatrix}0\\0\end{pmatrix},\begin{pmatrix}0\\-33\end{pmatrix},
                                                            \begin{pmatrix}-16\\30\end{pmatrix},\begin{pmatrix}44\\-33\end{pmatrix}
                                                            \right]$, see Figure \ref{fig:general_4_smallest}\\
         5 &                                      $=165$ & $\left[\begin{pmatrix}0\\0\end{pmatrix},
                                                           \begin{pmatrix}0\\-72\end{pmatrix},
                                                           \begin{pmatrix}-35\\12\end{pmatrix},
                                                           \begin{pmatrix}64\\-120\end{pmatrix},
                                                           \begin{pmatrix}-90\\-120\end{pmatrix}\right]$,
                                                           see Figure \ref{fig:general_5_smallest}\\
         6 &                                     $=1886$ & $\left[\begin{pmatrix}0\\0\end{pmatrix},
                                                           \begin{pmatrix}0\\-828\end{pmatrix},
                                                           \begin{pmatrix}-448\\-414\end{pmatrix},
                                                           \begin{pmatrix}-720\\132\end{pmatrix},
                                                           \begin{pmatrix}-1260\\-1023\end{pmatrix},
                                                           \begin{pmatrix}-1840\\-414\end{pmatrix}\right]$,\\
           &                                             & see Figure \ref{fig:general_6_smallest}\\
         7 &                                   $>599000$ & \\
      \hline
    \end{tabular}
    \caption{Best known constructions for maximal integral point sets over $\mathbb{Z}^2$ in general position.}
    \label{table_min_general_position}
  \end{center}
\end{table}

\section{Conclusion and outlook}
\label{sec_outlook}

\noindent
We have described several constructions for integral point sets over $\mathbb{Z}^2$ with given cardinality that fulfill some further properties. Although the maximality of the resulting integral point sets cannot be guaranteed so far, we conjecture them to be in many cases. We have described efficient algorithms for exhaustive generation of maximal integral point sets over $\mathbb{Z}^2$ and for testing the maximality of a given integral point set. Some exact values of minimum diameters for given cardinalities
could be obtained and several values are constructed as upper bounds and conjectured to be the exact values.

It remains a task to prove the maximality of point sets resulting from some of our constructions in general. Clearly similar problems could be considered in higher dimensions.

\bibliographystyle{amsplain}
\bibliography{mips}

\end{document}